\documentclass[12pt,leqno]{article}
\usepackage{amssymb}
\usepackage[mathscr]{eucal}
\usepackage{amsmath,amssymb,latexsym,theorem,bbm}
\usepackage{color}
\usepackage{appendix, url}

\newcommand{\CC}{\mathbb{C}}

\newcommand{\NN}{\mathbb{N}}

\newcommand{\RR}{\mathbb{R}}

\newcommand{\ZZ}{\mathbb{Z}}

\newcommand{\bA}{{\boldsymbol{A}}}

\newcommand{\bB}{{\boldsymbol{B}}}

\newcommand{\bc}{{\boldsymbol{c}}}

\newcommand{\bD}{{\boldsymbol{D}}}

\newcommand{\bI}{{\boldsymbol{I}}}

\newcommand{\bP}{{\boldsymbol{P}}}

\newcommand{\bQ}{{\boldsymbol{Q}}}

\newcommand{\bS}{{\boldsymbol{S}}}
\newcommand{\bs}{{\boldsymbol{s}}}

\newcommand{\bu}{{\boldsymbol{u}}}

\newcommand{\bV}{{\boldsymbol{V}}}

\newcommand{\bx}{{\boldsymbol{x}}}
\newcommand{\bX}{{\boldsymbol{X}}}
\newcommand{\by}{{\boldsymbol{y}}}
\newcommand{\bY}{{\boldsymbol{Y}}}

\newcommand{\bZ}{{\boldsymbol{Z}}}
\newcommand{\bU}{{\boldsymbol{U}}}

\newcommand{\Beta}{{\boldsymbol{\eta}}}

\newcommand{\btheta}{{\boldsymbol{\theta}}}

\newcommand{\bzeta}{{\boldsymbol{\zeta}}}
\newcommand{\bzero}{{\boldsymbol{0}}}

\newcommand{\cB}{{\mathcal B}}

\newcommand{\cD}{{\mathcal D}}
\newcommand{\cE}{{\mathcal E}}

\newcommand{\cG}{{\mathcal G}}
\newcommand{\cF}{{\mathcal F}}

\newcommand{\cN}{{\mathcal N}}

\newcommand{\dd}{\mathrm{d}}
\newcommand{\ee}{\mathrm{e}}
\newcommand{\ii}{\mathrm{i}}

\newcommand{\EE}{\operatorname{\mathbb{E}}}

\newcommand{\PP}{{\operatorname{\mathbb{P}}}}
\newcommand{\QQ}{{\operatorname{\mathbb{Q}}}}

\newcommand{\vare}{\varepsilon}

\renewcommand{\mid}{\,|\,}

\renewcommand{\leq}{\leqslant}
\renewcommand{\geq}{\geqslant}

\newcommand{\stoch}{\stackrel{\PP}{\longrightarrow}}
\newcommand{\stochG}{\stackrel{\PP_G}{\longrightarrow}}
\newcommand{\stochQ}{\stackrel{\QQ}{\longrightarrow}}
\newcommand{\stocheta}{\stackrel{\PP_{\{\exists\Beta^{-1}\}}}{\longrightarrow}}
\newcommand{\stochGeta}{\stackrel{\PP_{G\cap\{\exists\Beta^{-1}\}}}{\longrightarrow}}
\newcommand{\stochGketa}{\stackrel{\PP_{G_k\cap\{\exists\Beta^{-1}\}}}{\longrightarrow}}
\newcommand{\stochGetaHL}{\stackrel{\PP_{G\cap\{ \eta^2 > 0\}}}{\longrightarrow}}

\newcommand{\distrP}{\stackrel{\cD(\PP)}{\longrightarrow}}
\newcommand{\distrQ}{\stackrel{\cD(\QQ)}{\longrightarrow}}

\newcommand{\distre}{\stackrel{\cD}{=}}

\newcommand{\meanP}{\stackrel{L_1(\PP)}{\longrightarrow}}

\newcommand{\qmeanP}{\stackrel{L_2(\PP)}{\longrightarrow}}

\newcommand{\asP}{\stackrel{{\text{$\PP$-a.s.}}}{\longrightarrow}}
\newcommand{\asPG}{\stackrel{{\text{$\PP_G$-a.s.}}}{\longrightarrow}}

\newcommand{\aseP}{\stackrel{{\text{$\PP$-a.s.}}}{=}}

\newcommand{\bbone}{\mathbbm{1}}

\newcommand{\proofend}{\hfill\mbox{$\Box$}}

\numberwithin{equation}{section}

\theoremstyle{change} \theorembodyfont{\em}
\newtheorem{Lem}{Lemma.}[section]
\newtheorem{Thm}[Lem]{Theorem.}

\newtheorem{Cor}[Lem]{Corollary.}
\newtheorem{Def}[Lem]{Definition.}

\theorembodyfont{\rm}
\newtheorem{Rem}[Lem]{Remark.}

\begin{document}

\begin{center}
 {\bfseries\Large
   A multidimensional stable limit theorem}

\vspace*{3mm}

 {\sc\large
  M\'aty\'as $\text{Barczy}^{*}$,
  \ \framebox[1.1\width]{Gyula $\text{Pap}$}}

\end{center}

\vskip0.1cm

\noindent
 * ELKH-SZTE  Analysis and Applications Research Group,
   Bolyai Institute, University of Szeged,
   Aradi v\'ertan\'uk tere 1, H--6720 Szeged, Hungary.

\noindent e--mail: barczy@math.u-szeged.hu (M. Barczy).


\vskip0.5cm


\renewcommand{\thefootnote}{}
\vspace*{0.2cm}
\footnote{M\'aty\'as Barczy was supported by the project TKP2021-NVA-09.
Project no.\ TKP2021-NVA-09 has been implemented with the support
 provided by the Ministry of Innovation and Technology of Hungary from the National Research, Development and Innovation Fund,
 financed under the TKP2021-NVA funding scheme.}

\vspace*{-10mm}

\begin{abstract}
We establish multidimensional analogues of  one-dimensional stable limit theorems due to H\"ausler and
 Luschgy (2015) for so called explosive processes.
As special cases we present multidimensional stable limit theorems involving multidimensional normal-, Cauchy- and stable distributions as well.
\end{abstract}

\section{Introduction and main results}
\label{section_intro}

Stable convergence and mixing convergence have been frequently used in limit theorems in probability theory and statistics.
Historically the notion of mixing convergence was introduced first, and it can be traced back at least to R\'enyi \cite{Ren1},
 see also  R\'enyi \cite{Ren2} and \cite{RenRev}.
The more general concept of stable convergence is also due to R\'enyi \cite{Ren3}.
Stable convergence should not be mistaken for weak convergence to a stable distribution.
Recently, H\"ausler and Luschgy \cite{HauLus} have given an up to date and rigorous exposition of the mathematical theory
 of stable convergence, and they provided many applications in different areas to demonstrate the usefulness of this
 mode of convergence as well.
In many classical limit theorems, such as in the classical central limit theorem, not only convergence in distribution,
 but stable convergence holds as well, see, e.g., Examples 3.13 and 3.14 in H\"ausler and Luschgy \cite{HauLus}.
Stable convergence comes into play in the description of limit points of random sequences, in limit theorems with
 random indices, there is a version of the classical $\Delta$-method with stable convergence as well, see, e.g.,
 Chapter 4 in H\"ausler and Luschgy \cite{HauLus}.
Stable convergence has a central role in limit theorems for martingale difference arrays, and one can find its
 nice applications in describing the asymptotic behaviour of some estimators (such as conditional least squares estimator)
 of some parameters of autoregressive and moving average processes and supercritical Galton--Watson processes
 (for a detailed description, see Chapters 9 and 10 in H\"ausler and Luschgy \cite{HauLus}).
For a short survey on the role of stable convergence in limit theorems for semimartingales, see Podolskij and Vetter \cite{PodVet}.
In numerical probability, especially, in studying discretized processes, in approximation of stochastic integrals
 and stochastic differential equations, and in high frequency statistics, stable convergence also plays an essential role,
 see the recent books A\"it-Sahalia and Jacod \cite{AitJac} and Jacod and Protter \cite{JacPro}.
Very recently, Basse-O'Connor et al.\ \cite[part (i) of Theorem 2.1 and part (i) of Proposition 2.3]{BasHeiPod} have proved
 new limit theorems with stable convergence for some variational functionals of stationary increments L\'evy driven moving averages
 in the high frequency setting.

Recently, Crimaldi et al.\ \cite[Definition 3]{CriLetPra} have extended the notion of stable convergence:
 they have introduced the notion of stable convergence of random variables with respect to a so-called conditioning system towards a kernel,
 by replacing the single sub-$\sigma$-field appearing in the definition of (the original) stable
 convergence with a family of sub-$\sigma$-fields (called a conditioning system).
Then, as a generalization of the previously mentioned concept, Crimaldi et al.\ \cite[Definition 4]{CriLetPra}
 have introduced the notion of stable convergence of random variables in the strong sense with respect to a conditioning system,
 where not only the single sub-$\sigma$-field appearing in the definition of (the original) stable convergence
 is replaced by a conditioning system, but also the type of convergence for the conditional expectations
 with respect to the members of the conditional system in question is strengthened to convergence in $L_1$.
\ Moreover, as a further generalization, Crimaldi \cite[Definition 2.1]{Cri} have defined the notion of
 almost sure conditional convergence of random variables with respect to a conditional system towards a kernel.
If such a convergence holds, then the conditional expectations with respect to the members of the conditional system in question
 converge almost surely to a random variable.

Let \ $\ZZ_+$, \ $\NN$, \ $\RR$, \ $\RR_+$ \ and \ $\RR_{++}$ \ denote the set
 of non-negative integers, positive integers, real numbers, non-negative real
 numbers and positive real numbers, respectively.
The imaginary unit is denoted by \ $\ii$.
\ The Borel \ $\sigma$-algebra on \ $\RR^d$ \ is denoted by \ $\cB(\RR^d)$, \ where \ $d\in\NN$.
\ Further, let \ $\log^+(x):= \log(x)\bbone_{\{x\geq 1\}} + 0\cdot \bbone_{\{ 0\leq x < 1\}}$ \ for \ $x\in\RR_+$.
\ Convergence in a probability, in \ $L_1$, \ in \ $L_2$ \ and in distribution under a probability measure
 \ $\PP$ \ will be denoted by \ $\stoch$, \ $\meanP$, \ $\qmeanP$ \ and \ $\distrP$, \ respectively.
For an event \ $A$ \ with \ $\PP(A) > 0$,
 \ let \ $\PP_A(\cdot) := \PP(\cdot\mid A) = \PP(\cdot \cap A) / \PP(A)$ \ denote
 the conditional probability measure given \ $A$.
\ Let \ $\EE_\PP$ \ denote the expectation under a probability measure \ $\PP$.
\ Almost sure equality under a probability measure \ $\PP$ \ and equality in distribution will be denoted by \ $\aseP$
 \ and \ $\distre$, \ respectively.
Every random variable will be defined on a (suitable) probability space \ $(\Omega,\cF,\PP)$.
For a random variable \ $\xi:\Omega\to\RR^d$, \ the distribution of \ $\xi$ \ on \ $(\RR^d,\cB(\RR^d))$ \ is denoted by \ $\PP^\xi$.
\ The notions of stable and mixing convergence and some of their important properties used in the present paper are recalled in Appendix \ref{App_1}.

First, we will recall a one-dimensional stable limit theorem due to H\"ausler and Luschgy \cite[Theorem 8.2]{HauLus} for
 so called explosive processes.
The increments of these processes are in general not asymptotically negligible and do not satisfy
 the conditional Lindeberg condition, so they are not in the scope of stable martingale central limit theorems.
For such explosive processes, H\"ausler and Luschgy \cite{HauLus} developed the following limit theorem (Theorem \ref{THM_HL_8_2})
 which states stable (mixing) convergence of the appropriately scaled explosive process in question, and they successfully applied it for proving
 stable (mixing) convergence of conditional least squares estimator of the autoregressive parameter of supercritical
 autoregressive processes of order 1 (see H\"ausler and Luschgy \cite[Example 8.10 and Theorem 9.2]{HauLus})
 and that of Lotka-Nagaev estimator, conditional least squares estimator and Harris estimator of the
 offspring mean of supercritical Galton-Watson branching processes conditionally on non-extinction
 (see H\"ausler and Luschgy \cite[Corollaries 10.2, 10.4 and 10.6]{HauLus}).

\begin{Thm}[H\"ausler and Luschgy {\cite[Theorem 8.2]{HauLus}}]\label{THM_HL_8_2}
Let \ $(X_n)_{n\in\ZZ_+}$ \ and \ $(A_n)_{n\in\ZZ_+}$ \ be real-valued stochastic processes
 defined on a probability space \ $(\Omega,\cF,\PP)$ \ and adapted to a filtration \ $(\cF_n)_{n\in\ZZ_+}$.
\ Suppose that \ $A_n\geq 0$, \ $n\in\NN$, \ and that there exists \ $n_0\in\NN$ \ such that \ $A_n>0$ \ for each \ $n\geq n_0$.
\ Let \ $(a_n)_{n\in\NN}$ \ be a sequence in \ $(0,\infty)$ \ such that \ $a_n\to\infty$ \ as \ $n\to\infty$,
\ and let \ $G\in\cF_\infty:=\sigma(\bigcup_{n\in\ZZ_+} \cF_n)$ \ such that \ $\PP(G)>0$.
\ Assume that the following conditions are satisfied:
 \begin{enumerate}
  \item[\textup{(HLi)}]
    there exists a non-negative, $\cF_\infty$-measurable random variable \ $\eta:\Omega\to\RR$ \ such that
    \ $\PP(G\cap \{\eta^2 >0\})>0$ \ and
    \[
       \frac{A_n}{a_n^2}\stochG \eta^2 \qquad \text{as \ $n\to\infty$,}
    \]
  \item[\textup{(HLii)}]
   $(\frac{X_n}{a_n})_{n\in\NN}$ \ is stochastically bounded in \ $\PP_{G\cap \{\eta^2 >0\}}$-probability, i.e.,
   \[
     \lim_{K\to \infty} \sup_{n\in\NN}  \PP_{G\cap \{\eta^2 >0\}}\left( \frac{\vert X_n\vert}{a_n} > K\right) = 0,
   \]
  \item[\textup{(HLiii)}]
   there exists \ $p\in(1,\infty)$ \ such that
   \[
     \lim_{n\to\infty} \frac{a_{n-r}^2}{a_n^2} = \frac{1}{p^r} \qquad \text{for each \ $r\in\NN$,}
   \]
  \item[\textup{(HLiv)}]
     there exists a probability measure \ $\mu$ \ on \ $(\RR,\cB(\RR))$ \ with \ $\int_\RR \log^+(\vert x\vert)\,\mu(\dd x) < \infty$ \
     such that
     \[
      \EE_\PP\left( \exp\left\{ \ii t \frac{\Delta X_n}{A_n^{1/2}}\right\} \,\Big\vert\, \cF_{n-1}\right)
        \stochGetaHL \int_\RR \ee^{\ii tx}\,\dd \mu(x) \qquad \text{as \ $n\to\infty$}
     \]
     for all \ $t\in\RR$, \ where \ $\Delta X_n := X_n - X_{n-1}$, \ $n\in\NN$, \ and \ $\Delta X_0:=0$.
 \end{enumerate}
Then
 \begin{equation}\label{HL_BU}
  \frac{X_n}{A_n^{1/2}} \to \sum_{j=0}^\infty p^{-j/2} Z_j \qquad
  \text{$\cF_\infty$-mixing under \ $\PP_{G \cap \{\eta^2 > 0\}}$ \ as \ $n \to \infty$,}
 \end{equation}
 and
 \begin{equation}\label{HL_QU}
   \frac{X_n}{a_n} \to \eta \sum_{j=0}^\infty p^{-j/2} Z_j \qquad
   \text{$\cF_\infty$-stably under \ $\PP_{G \cap \{ \eta^2 >0 \}}$ \ as \ $n \to \infty$,}
 \end{equation}
 where \ $(Z_j)_{j\in\ZZ_+}$ \ denotes a \ $\PP$-independent and identically distributed sequence of real-valued random variables
 being \ $\PP$-independent of \ $\cF_\infty$ \ with \ $\PP(Z_0 \in B) = \mu(B)$ \ for all \ $B \in \cB(\RR)$.
\end{Thm}

\begin{Rem}
(i) The series \ $\sum_{j=0}^\infty p^{-j/2} Z_j = \sum_{j=0}^\infty (p^{1/2})^{-j} Z_j$ \ in \eqref{HL_BU} and \eqref{HL_QU} is absolutely convergent \ $\PP$-almost surely,
 since \ $p^{1/2}>1$, \ $\EE_\PP(\log^+(\vert Z_0\vert))<\infty$ \ (by condition (HLiv) of Theorem \ref{THM_HL_8_2}), and one can apply
 Lemma 8.1 in H\"ausler and Luschgy \cite{HauLus}.

(ii) We note that in condition (HLi) of Theorem \ref{THM_HL_8_2} the \ $\cF_\infty$-measurability of \ $\eta$ \ is supposed,
 but in condition (i) of Theorem 8.2 in H\"ausler and Luschgy \cite{HauLus} it is not supposed.
However, in the proof of Theorem 8.2 in H\"ausler and Luschgy \cite[page 148]{HauLus} it is written that
 the \ $\cF_\infty$-measurability of \ $\eta$ \ can be assumed without loss of generality.
Note also that if the probability space \ $(\Omega,\cF_\infty,\PP_G)$ \ is complete, then
 the \ $\cF_\infty$-measurability of \ $\eta$ \ follows itself from
 the convergence \ $\frac{A_n}{a_n^2}\stochG \eta^2$ \ as \ $n\to\infty$ \ involved in condition (HLi) of Theorem \ref{THM_HL_8_2}.
Indeed, then there exists a subsequence \ $(n_k)_{k\in\NN}$ \ such that
 \ $A_{n_k}/a_{n_k}^2 \asPG \eta^2$ \ as \ $k\to\infty$.
\ Since \ $A_{n_k}/a_{n_k}^2$ \ is \ $\cF_\infty$-measurable for each \ $k\in\NN$ \ and  \ $(\Omega,\cF_\infty,\PP_G)$ \ is complete,
 by a standard measure theoretical argument, we have \ $\eta^2$ \ is \ $\cF_\infty$-measurable.
The continuity of the square-root function together with \ $\eta\geq 0$ \ yield the \ $\cF_\infty$-measurability of \ $\eta$, \ as desired.

(iii)
The \ $\cF_\infty$-measurability of \ $\eta$ \ yields that \ $\eta$ \ and \ $Z_j$, \ $j\in\NN$, \ are \ $\PP$-independent
 in Theorem \ref{THM_HL_8_2}.
Further, we have \ $\PP_G(\eta>0)=\PP_G(\eta^2>0)>0$ \ and \ $\PP_{G\cap\{\eta^2>0\}}(\eta>0)=1$, where we used that \ $\eta$ \ is non-negative.
\proofend
\end{Rem}

By \ $\|\bx\|$ \ and \ $\|\bA\|$, \ we denote the Euclidean norm of a vector
 \ $\bx \in \RR^d$ \ and the induced matrix norm of a matrix
 \ $\bA \in \RR^{d\times d}$, \ respectively.
By \ $\langle \bx,\by\rangle$, \ we denote the Euclidean inner product of vectors \ $\bx,\by\in\RR^d$.
\ The null vector and the null matrix will be denoted by \ $\bzero$.
\ By \ $\varrho(\bA)$, \ we denote the spectral radius of \ $\bA \in \RR^{d\times d}$.
\ Moreover, \ $\bI_d \in \RR^{d\times d}$ \ denotes the $d\times d$ identity matrix, and if \ $\bA \in \RR^{d\times d}$ \ is symmetric and positive semidefinite,
 then \ $\bA^{1/2}$ \ denotes the unique symmetric, positive semidefinite square root of \ $\bA$.
\ If \ $\bV \in \RR^{d\times d}$ \ is symmetric and positive semidefinite, then
 \ $\cN_d(\bzero, \bV)$ \ denotes the $d$-dimensional normal distribution with mean vector \ $\bzero\in\RR^d$ \ and covariance matrix \ $\bV$.

In order to formulate our multidimensional stable limit theorems, we need the following result, which is a multidimensional generalization
 of Lemma 8.1 in H\"ausler and Luschgy \cite{HauLus}, and it is interesting on its own right.

\begin{Lem}\label{series_convergence}
Let \ $(\bZ_j)_{j\in\ZZ_+}$ \ be a \ $\PP$-independent and identically distributed sequence of \ $\RR^d$-valued random vectors.
Let \ $\bP \in \RR^{d\times d}$ \ be an invertible matrix with \ $\varrho(\bP) < 1$.
\ Then the following assertions are equivalent:
 \begin{enumerate}
  \item[\textup{(i)}]
   $\EE_\PP(\log^+(\|\bZ_0\|)) < \infty$.
  \item[\textup{(ii)}]
   $\sum_{j=0}^\infty \|\bP^j \bZ_j\| < \infty$ \ $\PP$-almost surely.
  \item[\textup{(iii)}]
   $\sum_{j=0}^\infty \bP^j \bZ_j$ \ converges \ $\PP$-almost surely in \ $\RR^d$.
  \item[\textup{(iv)}]
   $\bP^j \bZ_j \to \bzero$ \ as \ $j \to \infty$ \ $\PP$-almost surely.
 \end{enumerate}
\end{Lem}

The proof of Lemma \ref{series_convergence} and the proofs of all the forthcoming results can be found in Section \ref{Sec_Proofs}.
We note that from the proof of Lemma \ref{series_convergence} it turns out that for the implications \ (i) $\Rightarrow$ (ii) $\Rightarrow$  (iii) $\Rightarrow$ (iv),
  we do not need the invertibility of \ $\bP$, \ we only need it for \ (iv) $\Rightarrow$ (i).

For an \ $\RR^d$-valued stochastic process \ $(\bU_n)_{n\in\ZZ_+}$, \ the increments
 \ $\Delta \bU_n$, \ $n \in \ZZ_+$, \ are defined by \ $\Delta \bU_0 := \bzero$ \ and
 \ $\Delta \bU_n := \bU_n - \bU_{n-1}$ \ for \ $n \in \NN$.

Our main result is the following multidimensional analogue of Theorem 8.2 in H\"ausler and Luschgy \cite{HauLus} (see also Theorem \ref{THM_HL_8_2}).

\begin{Thm}\label{MSLTES}
Let \ $(\bU_n)_{n\in\ZZ_+}$ \ and \ $(\bB_n)_{n\in\ZZ_+}$ \ be \ $\RR^d$-valued and \ $\RR^{d\times d}$-valued stochastic processes,
 respectively, defined on a probability space \ $(\Omega,\cF,\PP)$ \ and adapted to a filtration \ $(\cF_n)_{n\in\ZZ_+}$.
\ Suppose that \ $\bB_n$ \ is invertible for sufficiently large \ $n \in \NN$.
\ Let \ $(\bQ_n)_{n\in\NN}$ \ be a sequence in \ $\RR^{d\times d}$ \ such that
 \ $\bQ_n \to \bzero$ \ as \ $n \to \infty$ \ and \ $\bQ_n$ \ is invertible for
 sufficiently large \ $n \in \NN$.
\ Let \ $G \in \cF_\infty := \sigma(\bigcup_{k=0}^\infty \cF_k)$ \ with \ $\PP(G) > 0$.
\ Assume that the following conditions are satisfied:
 \begin{enumerate}
  \item[\textup{(i)}]
   there exists an \ $\RR^{d\times d}$-valued, $\cF_\infty$-measurable random matrix \ $\Beta:\Omega\to\RR^{d\times d}$ \ such that
    \ $\PP(G \cap \{\exists\,\Beta^{-1}\}) > 0$ \ and
    \[
      \bQ_n \bB_n^{-1} \stochG \Beta \qquad \text{as \ $n \to \infty$,}
    \]
  \item[\textup{(ii)}]
   $(\bQ_n \bU_n)_{n\in\NN}$ \ is stochastically bounded in \ $\PP_{G\cap\{\exists\,\Beta^{-1}\}}$-probability, i.e.,
   \[
     \lim_{K\to\infty} \sup_{n\in\NN} \PP_{G\cap\{\exists\,\Beta^{-1}\}}(\|\bQ_n \bU_n\| > K) = 0,
   \]
  \item[\textup{(iii)}]
   there exists an invertible matrix \ $\bP \in \RR^{d\times d}$ \ with \ $\varrho(\bP) < 1$
    \ such that
    \[
      \bB_n \bB_{n-r}^{-1} \stochG \bP^r \qquad
      \text{as \ $n \to \infty$ \ for each \ $r \in \NN$,}
    \]
  \item[\textup{(iv)}]
   there exists a probability measure \ $\mu$ \ on \ $(\RR^d, \cB(\RR^d))$ \ with
    \ $\int_{\RR^d} \log^+(\|\bx\|) \, \mu(\dd\bx) < \infty$
    \ such that
    \[
      \EE_\PP\bigl(\ee^{\ii\langle\btheta,\bB_n\Delta\bU_n\rangle}
                     \mid \cF_{n-1}\bigr)
      \stochGeta \int_{\RR^d} \ee^{\ii\langle\btheta,\bx\rangle} \, \mu(\dd\bx)
      \qquad \text{as \ $n \to \infty$}
    \]
    for all \ $\btheta \in \RR^d$.
 \end{enumerate}
Then
 \begin{equation}\label{conv_BU}
  \bB_n \bU_n \to \sum_{j=0}^\infty \bP^j \bZ_j \qquad
  \text{$\cF_\infty$-mixing under \ $\PP_{G \cap \{\exists\,\Beta^{-1}\}}$ \ as \ $n \to \infty$,}
 \end{equation}
 and
 \begin{equation}\label{conv_QU}
   \bQ_n \bU_n \to \Beta \sum_{j=0}^\infty \bP^j \bZ_j \qquad
   \text{$\cF_\infty$-stably under \ $\PP_{G \cap \{\exists\,\Beta^{-1}\}}$ \ as \ $n \to \infty$,}
 \end{equation}
 where \ $(\bZ_j)_{j\in\ZZ_+}$ \ denotes a \ $\PP$-independent and identically distributed sequence of \ $\RR^d$-valued random vectors
 being \ $\PP$-independent of \ $\cF_\infty$ \ with \ $\PP(\bZ_0 \in B) = \mu(B)$ \ for all \ $B \in \cB(\RR^d)$.
\end{Thm}

\begin{Rem}\label{Rem_1}
(i) The series \ $\sum_{j=0}^\infty \bP^j \bZ_j$ \ in \eqref{conv_BU} and in \eqref{conv_QU} is absolutely convergent $\PP$-almost surely,
 since \ $\bP$ is invertible, \ $\varrho(\bP)<1$, \ $\EE_\PP(\log^+(\|\bZ_0\|))<\infty$ \ (by condition (iv) of Theorem \ref{MSLTES}),
 and one can apply Lemma \ref{series_convergence}.

(ii) The random variable \ $\Beta$ \ and the sequence \ $(\bZ_j)_{j\in\ZZ_+}$ \ are \ $\PP$-independent in Theorem \ref{MSLTES},
 since \ $\Beta$ \ is \ $\cF_\infty$-measurable and the sequence \ $(\bZ_j)_{j\in\ZZ_+}$ \ is \ $\PP$-independent of \ $\cF_\infty$.
Further, we have \ $\PP_G(\exists\,\Beta^{-1})>0$ \ and \ $\PP_{G\cap\{\exists\,\Beta^{-1}\}}(\exists\,\Beta^{-1})=1$.

(iii) The proof of Theorem \ref{MSLTES} (which can be found in Section \ref{Sec_Proofs}) follows
 the method of that of Theorem 8.2 in H\"ausler and Luschgy \cite{HauLus}. However, a natural question
 also occurs, namely, would it be possible to prove Theorem \ref{MSLTES} using the Cram\'er-Wold theorem
 for stable convergence (see, e.g., H\"ausler and Luschgy \cite[Corollary 3.19]{HauLus})?
We do not know the answer to this question.
The Cram\'er-Wold theorem for stable convergence states that, given $\RR^d$-valued random variables \ $\bX_n$, $n\in\NN$, \ and \ $\bX$,
 \ $\bX_n$ \ converges \ $\cG$-stably to \ $\bX$ \ as \ $n\to\infty$ \
 if and only if for all \ $\bu\in\RR^d$, \ the real-valued random variables
 \ $\langle \bu,\bX_n\rangle$ \ converges \ $\cG$-stably to the real-valued random variable
 \ $\langle \bu,\bX\rangle$ \ as \ $n\to\infty$ \ (where we used the setup given
 in Definition \ref{Def_HL_stable_conv}).
Here we only note that even in the proofs of multivariate central limit theorems with
 scaling matrices not converging to a fixed positive definite matrix,
 not only the Cram\'er-Wold theorem (for convergence in distribution) comes into play,
 but a key lemma originated to Bolthausan \cite{Bol} and its generalization due to Biscio et al.\ \cite[Lemma 3.2]{BisPoiWaa},
 for more details see Biscio et al.\ \cite{BisPoiWaa}.
 \proofend
\end{Rem}

In the next remark we reformulate condition (iii) of Theorem \ref{MSLTES} in the one-dimensional case.

\begin{Rem}
In case of \ $d=1$ \ (so not using boldface style in this case),
 if condition (i) of Theorem \ref{MSLTES} and \ $\PP(\exists\,\eta^{-1})=\PP(\eta\ne 0)=1$ \ hold, then condition (iii) of Theorem \ref{MSLTES} is
 equivalent to the following condition:
 \begin{align}\label{cond_iii_equiv}
 \text{there exists $P\in(-1,1)\setminus\{0\}$ such that $Q_nQ_{n-r}^{-1}\to P^r$ as $n\to\infty$ for each $r\in\NN$.}
 \end{align}
Indeed, if conditions (i) and (iii) of Theorem \ref{MSLTES} with \ $d=1$ \ and \ $\PP(\exists\,\eta^{-1})=1$ \ hold,
 then there exists \ $P\in(-1,1)\setminus\{0\}$ \ such that for each \ $r\in\NN$, \ we have
 \[
   Q_n Q_{n-r}^{-1} = Q_n B_n^{-1} B_n B_{n-r}^{-1}B_{n-r} Q_{n-r}^{-1} \stochG \eta P^r \eta^{-1} = P^r
     \qquad \text{as \ $n\to\infty$.}
 \]
Since \ $Q_n Q_{n-r}^{-1}$ \ is non-random, we have \eqref{cond_iii_equiv}.
Conversely, if condition (i) of Theorem \ref{MSLTES} with \ $d=1$, \ \ $\PP(\exists\,\eta^{-1})=\PP(\eta\ne 0)=1$, \ and \eqref{cond_iii_equiv} hold, then
 there exists \ $P\in(-1,1)\setminus\{0\}$ \ such that for each \ $r\in\NN$, \ we have
 \[
  B_nB_{n-r}^{-1} = B_n Q_n^{-1} Q_n Q_{n-r}^{-1} Q_{n-r} B_{n-r}^{-1}
    \stochG  \eta^{-1} P^r \eta = P^r
    \qquad \text{as \ $n\to\infty$,}
 \]
 i.e., condition (iii) of Theorem \ref{MSLTES} with \ $d=1$ \ holds.
Finally, note that, with the notation $a_n:=Q_n^{-1}$, condition \eqref{cond_iii_equiv} implies that
 for each \ $r\in\NN$ \ we have
 \[
  \frac{a_{n-r}^2}{a_n^2} = Q_n^2 Q_{n-r}^{-2} \to P^{2r} = \frac{1}{((P^2)^{-1})^r}
   \qquad \text{as \ $n\to\infty$,}
 \]
 which is nothing else but condition (iii) of Theorem 8.2 in H\"ausler and Luschgy \cite{HauLus}
 (see also condition (HLiii) of Theorem \ref{THM_HL_8_2}) with $p:=(P^2)^{-1}\in(1,\infty)$.
In Remark \ref{Rem_comparison}, we give a more detailed comparison of Theorem 8.2 in H\"ausler and Luschgy \cite{HauLus}
 (see also Theorem \ref{THM_HL_8_2}) and Theorem \ref{MSLTES}.
\proofend
\end{Rem}

In the next remark we investigate the connection between Theorem 8.2 in H\"ausler and Luschgy \cite{HauLus} (see also Theorem \ref{THM_HL_8_2})
 and Theorem \ref{MSLTES}.

\begin{Rem}\label{Rem_comparison}
Theorem \ref{MSLTES} gives back Theorem 8.2 in H\"ausler and Luschgy \cite{HauLus} (see also Theorem \ref{THM_HL_8_2})
 provided that \ $\PP(\eta>0)=1$ \ in condition (i) of Theorem 8.2 in H\"ausler and Luschgy \cite{HauLus}.
Indeed, let \ $(X_n)_{n\in\ZZ_+}$ \ and \ $(A_n)_{n\in\ZZ_+}$ \ be real-valued stochastic processes defined on a probability
 space \ $(\Omega,\cF,\PP)$ \ and adapted to a filtration \ $(\cF_n)_{n\in\ZZ_+}$. \
Suppose that \ $A_n \geq 0$, \ $n\in\NN$, \ and that there exists \ $n_0\in\NN$ \ such that \ $A_n>0$ \ for each \ $n\geq n_0$.
 \ Let \ $(a_n)_{n\in\NN}$ \ be a sequence in \ $(0,\infty)$ \ such that \ $a_n\to\infty$ \ as \ $n\to\infty$, \ and let \ $G\in \cF_\infty$ \ with
 \ $\PP(G)>0$ \ such that the conditions (HLi) together with \ $\PP(\eta > 0)=1$, \ (HLii), (HLiii) and (HLiv) of Theorem \ref{THM_HL_8_2} hold.
Note that in this case \ $\PP_{G\cap\{\eta^2>0\}} = \PP_G$, \ since \ $\PP(\eta>0)=1$ \ implies that \ $\PP(\eta^2>0)=1$.
\ In Theorem \ref{MSLTES}, let us make the following choices
 \ $\bU_n := X_n$, \ $n\in\ZZ_+$, \ $\bB_n := A_n^{-1/2}$, \ $n\geq n_0$, \ $\bQ_n := a_n^{-1}$, \ $n\in\NN$, \ and \ $\bP:=p^{-1/2}$,
 \ where \ $p\in(1,\infty)$ \ is given in (HLiii) of Theorem \ref{THM_HL_8_2}.
Then (HLi) of Theorem \ref{THM_HL_8_2}, the non-negativity of \ $\eta$ \ and the continuity of the square-root function yield that
 \ $\bQ_n \bB_n^{-1} = \frac{A_n^{1/2}}{a_n} \stochG \eta$ \ as \ $n\to\infty$,
 \ i.e., condition (i) of Theorem \ref{MSLTES} is satisfied.
Further, (HLi) of Theorem \ref{THM_HL_8_2} together with \ $\PP(\eta > 0)=1$, \ (HLiii) of Theorem \ref{THM_HL_8_2}
 and the continuity of the square-root function imply that for each \ $r\in\NN$, \ we have
 \[
    \bB_n \bB_{n-r}^{-1} = \frac{A_{n-r}^{1/2}}{A_n^{1/2}}
                     = \frac{A_{n-r}^{1/2}}{a_{n-r}} \frac{a_n}{A_n^{1/2}} \frac{a_{n-r}}{a_n}
                     \stochG \eta \cdot\frac{1}{\eta}\cdot \frac{1}{p^{r/2}}
                     = \bP^r
                     \qquad \text{as \ $n\to\infty$,}
 \]
 i.e., condition (iii) of Theorem \ref{MSLTES} holds.
Conditions (HLii) and (HLiv) of Theorem \ref{THM_HL_8_2} readily yield conditions (ii) and (iv) of Theorem \ref{MSLTES}, respectively.
So we can apply Theorem \ref{MSLTES} and we have \eqref{HL_BU} and \eqref{HL_QU}, as desired.
\proofend
\end{Rem}

Next, we present a multidimensional stable central limit theorem, which is a multidimensional counterpart
 of Corollary 8.5 in H\"ausler and Luschgy \cite{HauLus}.

\begin{Cor}\label{Cor_Gauss}
Let us assume that the conditions of Theorem \ref{MSLTES} hold with \ $\mu:=\PP^{\,\cN_d(\bzero,\bD)}$, \ where
 \ $\PP^{\,\cN_d(\bzero,\bD)}$ \ denotes the distribution of a \ $d$-dimensional normally distributed random variable
 with mean vector \ $\bzero\in\RR^d$ \ and covariance matrix \ $\bD\in\RR^{d\times d}$.
\ Then
 \begin{equation}\label{conv_BU_Gauss}
  \bB_n \bU_n \to \bZ \qquad
  \text{$\cF_\infty$-mixing under \ $\PP_{G \cap \{\exists\,\Beta^{-1}\}}$ \ as \ $n \to \infty$,}
 \end{equation}
 and
 \begin{equation}\label{conv_QU_Gauss}
   \bQ_n \bU_n \to \Beta \bZ \qquad
   \text{$\cF_\infty$-stably under \ $\PP_{G \cap \{\exists\,\Beta^{-1}\}}$ \ as \ $n \to \infty$,}
 \end{equation}
 where \ $\bZ$ \ denotes a \ $d$-dimensional normally distributed random vector with mean vector
 \ $\bzero\in\RR^d$ \ and covariance matrix \ $\sum_{j=0}^\infty \bP^j \bD (\bP^j)^\top$,
 and \ $\bZ$ \ is \ $\PP$-independent of \ $\cF_\infty$.
\end{Cor}

In Corollary \ref{Cor_Gauss}, \ $\Beta$ \ and \ $\bZ$ \ are \ $\PP$-independent, since \ $\Beta$ \ is \ $\cF_\infty$-measurable
 (supposed in condition (i) of Theorem \ref{MSLTES}).

Next, we will formulate a corollary of Theorem \ref{MSLTES} involving multidimensional stable distributions,
 in particular, a multidimensional Cauchy distribution.
For this, first we recall the notion of a multidimensional stable distribution.
A \ $d$-dimensional random variable \ $\bzeta:=(\zeta_1,\ldots,\zeta_d)$ \ is said to be stable if for any \ $a_1,a_2\in\RR_{++}$ \
 there exist \ $b\in\RR_{++}$ \ and \ $\bc\in\RR^d$ \ such that
 \begin{align}\label{def_multi_stable}
     a_1\bzeta^{(1)} + a_2\bzeta^{(2)} \distre b\bzeta + \bc,
 \end{align}
 where \ $\bzeta^{(1)}$ \ and  \ $\bzeta^{(2)}$ \ are independent copies of \ $\bzeta$.
\ It is known that \ $\bzeta$ \ is stable if and only if there exists \ $\alpha\in(0,2]$ \ such that for each \ $n\geq 2$, \ $n\in\NN$ \ there exists \ $\bc_n\in\RR^d$ \
 satisfying \ $\bzeta^{(1)} + \cdots +\bzeta^{(n)} \distre n^{\frac{1}{\alpha}}\bzeta + \bc_n$, \ where \ $\bzeta^{(1)},\bzeta^{(2)},\ldots,\bzeta^{(n)}$ \
 are independent copies of \ $\bzeta$.
\ The index \ $\alpha$ \ is called the index of stability or the characteristic exponent of \ $\bzeta$.
\ In what follows, let \ $S_{d-1}:=\{  \bx\in\RR^d : \Vert \bx\Vert=1\}$ \ be the unit surface in \ $\RR^d$.
\ We say that \ $\bzeta$ \ is symmetric stable if it is stable and \ $\PP(\bzeta \in A) = \PP(-\bzeta \in A)$ \ for all \ $A\in\cB(\RR^d)$.
\ It known that a \ $d$-dimensional random variable \ $\bzeta$ \ is symmetric $\alpha$-stable with index
 \ $\alpha\in(0,2)$  \ if and only if there exists a unique symmetric finite measure \ $\Pi$ \ on \ $(S_{d-1},\cB(S_{d-1}))$
 \ (where the property symmetric means that \ $\Pi(A) = \Pi(-A)$ \ for any \ $A\in\cB(S_{d-1})$) \ such that
 \[
   \EE_\PP\big( \exp(\ii \langle \btheta, \bzeta \rangle) \big)
      = \exp\left\{  -\int_{S_{d-1}} \big\vert \langle \btheta, \bx \rangle  \big\vert^\alpha\,\Pi(\dd \bx) \right\},
         \qquad \btheta\in\RR^d,
 \]
 see, e.g., Sato \cite[Theorem 14.13]{Sat}.
The measure \ $\Pi$ \ is called the spectral measure of \ $\bzeta$.
We say that a \ $d$-dimensional random variable \ $\bzeta$ \ has a \ $d$-dimensional Cauchy distribution
 with parameter \ $(\bzero,\bI_d)$, \ if its density function takes the form
 \[
 f_\bzeta(\bx) = \frac{\Gamma\left(\frac{1+d}{2}\right)}{\pi^{\frac{1+d}{2}}}
                                  \Big(1+\Vert \bx\Vert^2\Big)^{-\frac{1+d}{2}},
                                  \qquad \bx\in\RR^d,
 \]
 see, e.g., Kotz and Nadarajah \cite[Section 2.2, page 41]{KotNad} or Sato \cite[Example 2.12]{Sat}.
It is known that if \ $\bzeta$ \ has a \ $d$-dimensional Cauchy distribution with parameter \ $(\bzero,\bI_d)$,
 \ then the characteristic function of \ $\bzeta$ \ takes the form
 \ $\EE_\PP( \ee^{\ii \langle \btheta, \bzeta \rangle} )  = \ee^{-\Vert \btheta \Vert}$, \ $\btheta\in\RR^d$,
 \ and \ $\bzeta$ \ is symmetric $1$-stable, see, e.g., Sato \cite[Theorem 14.14]{Sat}.

\begin{Cor}\label{Cor_Stable}
Let us assume that the conditions of Theorem \ref{MSLTES} hold with \ $\mu:=\PP^{\bzeta}$, \ where
 \ $\bzeta$ \ is a \ $d$-dimensional symmetric \ $\alpha$-stable random variable with characteristic exponent
 \ $\alpha\in(0,2)$ \ and spectral measure \ $\Pi$.
\ Then
 \begin{equation}\label{conv_BU_Stable}
  \bB_n \bU_n \to \bZ \qquad
  \text{$\cF_\infty$-mixing under \ $\PP_{G \cap \{\exists\,\Beta^{-1}\}}$ \ as \ $n \to \infty$,}
 \end{equation}
 and
 \begin{equation}\label{conv_QU_Stable}
   \bQ_n \bU_n \to \Beta \bZ \qquad
   \text{$\cF_\infty$-stably under \ $\PP_{G \cap \{\exists\,\Beta^{-1}\}}$ \ as \ $n \to \infty$,}
 \end{equation}
 where \ $\bZ$ \ denotes a \ $d$-dimensional random vector \ $\PP$-independent of \ $\cF_\infty$ \ with a
 characteristic function
 \begin{align}\label{help_Stable_limit_char}
   \EE_\PP(\ee^{\ii\langle \btheta,\bZ \rangle})
     = \exp\left\{  - \int_{S_{d-1}}\sum_{j=0}^\infty \big\vert \langle (\bP^j)^\top \btheta,\bx \rangle \big\vert^\alpha \, \Pi(\dd \bx)\right\},
       \qquad \btheta\in\RR^d.
 \end{align}
In particular, if \ $\bzeta$ \ has a $d$-dimensional Cauchy distribution with parameter \ $(\bzero,\bI_d)$, \
 then \ $\bZ$ \ has a characteristic function
 \begin{align}\label{help_Cauchy_limit_char}
   \EE_\PP(\ee^{\ii\langle \btheta,\bZ \rangle})
     = \exp\left\{  - \sum_{j=0}^\infty \Vert (\bP^j)^\top \btheta \Vert \right\},
       \qquad \btheta\in\RR^d.
 \end{align}
\end{Cor}

In Corollary \ref{Cor_Stable}, \ $\Beta$ \ and \ $\bZ$ \ are \ $\PP$-independent, since \ $\Beta$ \ is \ $\cF_\infty$-measurable.
Corollary \ref{Cor_Stable} in the special case when $\bzeta$ has a $d$-dimensional Cauchy distribution with parameter $(\bzero,\bI_d)$
 can be considered as a multidimensional counterpart of Exercise 8.1 in H\"ausler and Luschgy \cite{HauLus}.

Finally, we formulate a slight generalization of Theorem \ref{MSLTES} in case of \ $G=\Omega$, \ by weakening its condition (iv) a little bit.
This generalization can be considered as a multidimensional analogue of Corollary 8.8 in H\"ausler and Luschgy \cite{HauLus}.

\begin{Cor}\label{MSLTES_discrete_S}
Let us suppose that the conditions of Theorem \ref{MSLTES} are satisfied with \ $G:=\Omega$ \ except its condition (iv) which is replaced by
 \begin{align*}
 &\text{(iv') \quad there exists a probability measure \ $\mu$ \ on \ $(\RR^d,\cB(\RR^d))$ \ with
        $\int_{\RR^d} \log^+(\Vert\bx\Vert) \,\mu(\dd\bx)<\infty$,}\\[-2mm]
 &\phantom{(iv') \quad}
   \text{and an \ $\cF_\infty$-measurable, \ $\RR^{d\times d}$-valued discrete random variable \ $\bS$ \ such that}\\
 &\phantom{(iv') \quad}
    \EE_\PP\bigl(\ee^{\ii\langle\btheta,\bB_n\Delta\bU_n\rangle}
                     \mid \cF_{n-1}\bigr)
      \stocheta \int_{\RR^d} \ee^{\ii\langle\btheta,\bS\bx\rangle} \, \mu(\dd\bx)
       \; \text{\ as \ $n \to \infty$ \ for all \ $\btheta \in \RR^d$.}
 \end{align*}
Then
 \begin{equation}\label{conv_BU_discrete_S}
  \bB_n \bU_n \to \sum_{j=0}^\infty \bP^j \bS\bZ_j \qquad
  \text{$\cF_\infty$-stably under \ $\PP_{\{\exists\,\Beta^{-1}\}}$ \ as \ $n \to \infty$,}
 \end{equation}
 and
 \begin{equation}\label{conv_QU_discrete_S}
   \bQ_n \bU_n \to \Beta \sum_{j=0}^\infty \bP^j \bS\bZ_j \qquad
   \text{$\cF_\infty$-stably under \ $\PP_{\{\exists\,\Beta^{-1}\}}$ \ as \ $n \to \infty$,}
 \end{equation}
 where \ $(\bZ_j)_{j\in\ZZ_+}$ \ denotes a \ $\PP$-independent and identically distributed sequence of \ $\RR^d$-valued random vectors
 \ $\PP$-independent of \ $\cF_\infty$ \ with \ $\PP(\bZ_0 \in B) = \mu(B)$ \ for all \ $B \in \cB(\RR^d)$.
\end{Cor}

In Corollary \ref{MSLTES_discrete_S}, \ $\Beta$ \ and \ $(\bZ_j)_{j\in\ZZ_+}$ \ are \ $\PP$-independent (see part (ii) of Remark \ref{Rem_1}).
For an application of Corollary \ref{MSLTES_discrete_S} with $d=1$, see the proof of
 Theorem 9.1 in H\"ausler and Luschgy \cite{HauLus}, where the authors prove stable convergence
 of conditional least squares estimator of the autoregressive parameter of supercritical
 autoregressive processes of order 1.

Finally, we note that in a companion paper Barczy and Pap \cite{BarPap}, we use our main result Theorem \ref{MSLTES} for studying the asymptotic behaviour of least squares
 estimator of the autoregressive parameters of some supercritical Gaussian autoregressive processes of order 2 using random scaling.
In another companion paper Barczy \cite{Bar}, we also use Theorem \ref{MSLTES} for proving stable convergence
 of conditional least squares estimators of drift parameters for supercritical continuous state and continuous time
 branching processes with immigration based on discrete time observations.

\section{Proofs}\label{Sec_Proofs}

\noindent
\textbf{Proof of Lemma \ref{series_convergence}.}
(i) $\Rightarrow$ (ii).
\ We have \ $\varrho(\bP) = \lim_{k\to\infty} \|\bP^k\|^{1/k}$ \ by the Gelfand formula, see,
 e.g., Horn and Johnson \cite[Corollary 5.6.14]{HorJoh}.
Hence there exists \ $k_0 \in \NN$ \ such that
 \begin{equation}\label{Gelfand}
  \|\bP^k\|^{1/k} \leq \varrho(\bP) + \frac{1-\varrho(\bP)}{2} = \frac{1+\varrho(\bP)}{2} < 1
  \qquad \text{for each \ $k \geq k_0$,}
 \end{equation}
 since \ $\varrho(\bP) < 1$.
\ Choose \ $c \in \bigl(1, \frac{2}{1+\varrho(\bP)}\bigr)$.
\ Then (i) implies
 \begin{align*}
   \sum_{j=k_0}^\infty \PP(\|\bZ_j\| > c^j)
  & = \sum_{j=k_0}^\infty \PP(\|\bZ_0\| > c^j)
   = \sum_{j=k_0}^\infty \PP(\log^+(\|\bZ_0\|) > j \log^+(c)) \\
  & =  \sum_{j=k_0}^\infty \PP\left( \frac{\log^+(\|\bZ_0\|)}{\log(c)} > j \right)
   < \infty ,
 \end{align*}
 where we used that \ $\log^+(c)=\log(c)>0$ \ and \ $\sum_{n=1}^\infty \PP(\xi\geq n)\leq \EE_\PP(\xi)$ \ for any non-negative random variable \ $\xi$.
\ By the Borel--Cantelli lemma,
 \[
   \PP\big(\limsup_{j\to\infty} \{\|\bZ_j\| > c^j\}\big) = 0,
   \qquad \text{and hence} \qquad
   \PP\big(\liminf_{j\to\infty} \{\|\bZ_j\| \leq c^j\}\big) = 1 ,
 \]
 i.e., for \ $\PP$-a.a.\ $\omega\in\Omega$, \ there exists \ $j_0(\omega)\in\NN$ \ such that \ $\|\bZ_j(\omega)\| \leq c^j$ \ for each \ $j\geq j_0(\omega)$.
\ Consequently, for \ $\PP$-a.a.\ $\omega\in\Omega$, \ we have
 \[
  \sum_{j=k_0 \vee j_0(\omega)}^\infty \|\bP^j \bZ_j(\omega)\|
    \leq \sum_{j=k_0 \vee j_0(\omega) }^\infty \|\bP^j\| \cdot \|\bZ_j(\omega)\|
    \leq \sum_{j=k_0}^\infty \left(\frac{1+\varrho(\bP)}{2}\right)^j c^j < \infty,
 \]
 since \ $\frac{1+\varrho(\bP)}{2}c\in(0,1)$.
\ It yields (ii).

The implications (ii) $\Rightarrow$ (iii) and (iii) $\Rightarrow$ (iv) are obvious.

(iv) $\Rightarrow$ (i).
\ We have \ $\PP(\limsup_{j\to\infty} \{\|\bP^j \bZ_j\| > 1\}) = 0$, \ and hence, by the Borel--Cantelli lemma
 and the independence of \ $(\bZ_j)_{j\in\ZZ_+}$, \ we get
 \[
   \sum_{j=0}^\infty \PP(\|\bP^j \bZ_j\| > 1) < \infty .
 \]
Using that the determinant of \ $\bP$ \ coincides with the product of its eigenvalues,
 the invertibility of \ $\bP$ \ implies that \ $\bP$ \ does not have an eigenvalue \ $0$, \
 and, in particular, we get \ $\varrho(\bP) > 0$.
\ The eigenvalues of \ $\bP^{-1}$ \ are the reciprocals of the eigenvalues of \ $\bP$, \
 hence \ $\varrho(\bP^{-1}) \geq \frac{1}{\varrho(\bP)}$, \ implying \ $\|\bP^{-1}\| \geq \varrho(\bP^{-1}) \geq \frac{1}{\varrho(\bP)} > 1$.
\ Thus for each \ $j \in \ZZ_+$, \ we have
 \begin{align*}
  \PP(\|\bP^j \bZ_j\| > 1)
  & = \PP(\|\bP^{-1}\|^j \|\bP^j \bZ_0\| > \|\bP^{-1}\|^j)
    \geq \PP(\|(\bP^{-1})^j \bP^j \bZ_0\| > \|\bP^{-1}\|^j) \\
  & = \PP(\|\bZ_0\| > \|\bP^{-1}\|^j)
    = \PP(\log^+(\|\bZ_0\|) > j \log^+(\|\bP^{-1}\|)) .
 \end{align*}
Consequently,
 \ $\sum_{j=0}^\infty \PP(\log^+(\|\bZ_0\|) > j \log(\|\bP^{-1}\|)) < \infty$, \ yielding
 \[
  \EE_\PP\left( \frac{\log^+(\|\bZ_0\|)}{\log(\|\bP^{-1}\|)} \right)<\infty
 \]
 and hence (i), where we used that \ $\log^+(\|\bP^{-1}\|) = \log(\|\bP^{-1}\|) >0$ \
 and \ $\EE_\PP(\xi) \leq 1 + \sum_{n=1}^\infty \PP(\xi>n)$ \ for any non-negative random variable \ $\xi$.
\proofend

\medskip

\noindent
\textbf{Proof of Theorem \ref{MSLTES}.}\\
{\sl Step 1:} Let \ $\QQ := \PP_{G\cap\{\exists\,\Beta^{-1}\}}$, \ and for each \ $n \in \ZZ_+$, \ put
 \[
   L_n := \frac{\PP(G \cap \{\exists\, \Beta^{-1}\} \mid \cF_n)}{\PP(G \cap \{\exists\, \Beta^{-1}\})} .
 \]
Then \ $\QQ$ \ is absolutely continuous with respect to \ $\PP$ \ and \ $\PP_G$ \ as well, and, for each \ $n\in\NN$,
 \ $L_n$ \ is a well-defined and \ $\cF_n$-measurable random variable, since \ $\PP(G\cap\{\exists\,\Beta^{-1}\}) > 0$.
\ Note that \ $(L_n)_{n\in\ZZ_+}$ \ is the density process of \ $\QQ$ \ with respect to
 \ $\PP$, \ that is, \ $L_n = \frac{\dd\,\QQ|_{\cF_n}}{\dd\,\PP|_{\cF_n}}$ \ for every
 \ $n \in \ZZ_+$, \ where \ $\QQ|_{\cF_n}$ \ and \ $\PP|_{\cF_n}$ \ denote the restriction of \ $\QQ$ \ and
  \ $\PP$ \ onto \ $(\Omega,\cF_n)$, \ respectively.
Indeed, for all \ $A \in \cF_n$, \ we have
 \begin{align*}
  \QQ|_{\cF_n}(A) = \QQ(A) = \frac{\PP(A \cap G \cap \{\exists\,\Beta^{-1}\})}{\PP(G \cap \{\exists\,\Beta^{-1}\}) } ,
 \end{align*}
 and, by the definition of conditional expectation with respect to the \ $\sigma$-algebra \ $\cF_n$,
 \begin{align*}
  \int_A L_n(\omega)\,\PP|_{\cF_n}(\dd\omega)
  & = \int_A \frac{\PP(G \cap \{\exists\, \Beta^{-1}\} \mid \cF_n)}{\PP(G \cap \{\exists\,\Beta^{-1}\})}(\omega) \,\PP|_{\cF_n}(\dd\omega) \\
  &= \frac{1}{\PP(G \cap \{\exists\,\Beta^{-1}\})} \int_A (\EE_\PP(\bbone_{G\cap\{\exists\,\Beta^{-1}\}} \mid \cF_n))(\omega) \, \PP(\dd\omega)\\
  & = \frac{1}{\PP(G \cap \{\exists\,\Beta^{-1}\})} \int_A \bbone_{G\cap\{\exists\,\Beta^{-1}\}}(\omega)\,\PP(\dd\omega)
    =  \frac{\PP(A\cap G \cap \{\exists\,\Beta^{-1}\})}{\PP(G \cap \{\exists\,\Beta^{-1}\})} ,
 \end{align*}
 yielding that \ $\QQ|_{\cF_n}(A) = \int_A L_n(\omega) \,\PP|_{\cF_n}(\dd\omega)$, \ $A\in\cF_n$, \ as desired.
Then, by L\'evy's upwards theorem (see, e.g., Theorem \ref{Thm_Levy01}), we get
 \begin{align}\label{MCT1}
  &L_n \meanP
          \frac{\EE_\PP( \bbone_{ G \cap \{\exists\,\Beta^{-1}\} }\mid \cF_\infty)}{\PP(G \cap \{\exists\,\Beta^{-1}\})}
          =\frac{\bbone_{G\cap\{\exists\,\Beta^{-1}\}}}{\PP(G \cap \{\exists\,\Beta^{-1}\})} = \frac{\dd\QQ}{\dd\PP}
  \qquad \text{as \ $n \to \infty$,}\\\label{MCT2}
  &  L_n \asP  \frac{\EE_\PP( \bbone_{ G \cap \{\exists\,\Beta^{-1}\} }\mid \cF_\infty)}{\PP(G \cap \{\exists\,\Beta^{-1}\})}
                = \frac{\bbone_{G\cap\{\exists\,\Beta^{-1}\}}}{\PP(G \cap \{\exists\,\Beta^{-1}\})} = \frac{\dd\QQ}{\dd\PP}\qquad  \text{as \ $n \to \infty$,}
 \end{align}
 where the second equality in \eqref{MCT1} (and in \eqref{MCT2}) holds, since for all \ $A\in\cF$,
  \[
    \QQ(A)= \PP_{G \cap \{\exists\,\Beta^{-1}\}}(A) = \frac{\PP(A\cap G \cap \{\exists\,\Beta^{-1}\})}
                                                         {\PP(G \cap \{\exists\,\Beta^{-1}\})},
  \]
  and
  \[
   \int_A \frac{\bbone_{G \cap \{\exists\,\Beta^{-1}\}}(\omega)}{\PP(G \cap \{\exists\,\Beta^{-1}\})} \,\PP(\dd\omega)
      = \frac{\PP(A\cap G \cap \{\exists\,\Beta^{-1}\})}
               {\PP(G \cap \{\exists\,\Beta^{-1}\})}.
  \]

Next, we check that \ $\bZ_j$, $j\in\ZZ_+$, \ and \ $\cF_\infty$ \ are independent under \ $\QQ$ \ as well.
Indeed, since \ $\bZ_j$, $j\in\ZZ_+$, \ and \ $\cF_\infty$ \ are independent under \ $\PP$ \ (by assumption)
 and \ $G \cap \{\exists\,\Beta^{-1}\} \in \cF_\infty$
 \ (since \ $G\in \cF_\infty$ \ and \ $\Beta$ \ is \ $\cF_\infty$-measurable),
 we have for each \ $m\in\NN$, \ $B_0,B_1,\ldots,B_m\in\cB(\RR^d)$ \ and \ $A\in\cF_\infty$,
 \begin{align*}
   \QQ&(\{\bZ_0\in B_0\}\cap\{\bZ_1\in B_1\}\cap\cdots \cap\{\bZ_m\in B_m\}\cap A) \\
   & = \frac{\PP(\{\bZ_0\in B_0\}\cap\{\bZ_1\in B_1\}\cap\cdots \cap\{\bZ_m\in B_m\}\cap A \cap G\cap \{\exists\,\Beta^{-1}\})}
              {\PP(G\cap \{\exists\,\Beta^{-1}\})}\\
   & = \frac{\PP(\{\bZ_0\in B_0\}\cap\{\bZ_1\in B_1\}\cap\cdots \cap\{\bZ_m\in B_m\}) \PP( A \cap G\cap \{\exists\,\Beta^{-1}\})}
              {\PP(G\cap \{\exists\,\Beta^{-1}\})}\\
   &= \frac{\PP(\bZ_0\in B_0)\PP(\bZ_1\in B_1)\cdots \PP(\bZ_m\in B_m) \PP( A \cap G\cap \{\exists\,\Beta^{-1}\})}
              {\PP(G\cap \{\exists\,\Beta^{-1}\})},
 \end{align*}
 and
 \begin{align*}
   \QQ&(\bZ_0\in B_0)\QQ(\bZ_1\in B_1)\cdots \QQ(\bZ_m\in B_m)\QQ(A) \\
      & = \frac{\PP(\{\bZ_0\in B_0\}\cap G\cap \{\exists\,\Beta^{-1}\})}{\PP(G\cap \{\exists\,\Beta^{-1}\})}
           \cdots \frac{\PP(\{\bZ_m\in B_m\}\cap G\cap \{\exists\,\Beta^{-1}\})}{\PP(G\cap \{\exists\,\Beta^{-1}\})}
           \cdot \frac{\PP(A\cap G\cap \{\exists\,\Beta^{-1}\})}{\PP(G\cap \{\exists\,\Beta^{-1}\})} \\
      & = \PP(\{\bZ_0\in B_0\}) \cdots \PP(\{\bZ_m\in B_m\}) \frac{\PP(A\cap G\cap \{\exists\,\Beta^{-1}\})}{\PP(G\cap \{\exists\,\Beta^{-1}\})},
 \end{align*}
 where we used that
 \[
    \PP(\{\bZ_j\in B_j\}\cap G\cap \{\exists\,\Beta^{-1}\})
       = \PP(\bZ_j\in B_j) \PP(G\cap \{\exists\,\Beta^{-1}\}),\qquad j\in\{0,1,\ldots,m\}.
 \]
It yields that
 \begin{align*}
  &\QQ(\{\bZ_0\in B_0\}\cap\{\bZ_1\in B_1\}\cap\cdots \cap\{\bZ_m\in B_m\}\cap A) \\
  &\qquad = \QQ(\bZ_0\in B_0)\QQ(\bZ_1\in B_1)\cdots \QQ(\bZ_m\in B_m)\QQ(A),
 \end{align*}
 as desired.

For each \ $\btheta \in \RR^d$, \ let us introduce the notation
 \begin{align}\label{help1}
   \varphi_\mu(\btheta) := \int_{\RR^d} \ee^{\ii\langle\btheta,\bx\rangle} \, \mu(\dd\bx)
   = \EE_\PP(\ee^{\ii\langle\btheta,\bZ_0\rangle})
   = \EE_\QQ(\ee^{\ii\langle\btheta,\bZ_0\rangle}) ,
 \end{align}
 since the distributions of \ $\bZ_0$ \ under \ $\PP$ \ and \ $\QQ$ \ coincide.
\ Indeed, by the independence of \ $\cF_\infty$ \ and \ $\bZ_0$ \ under \ $\PP$, \ for all \ $B \in \cB(\RR^d)$, \ we have
 \[
   \QQ(\bZ_0 \in B) = \frac{\PP(\{\bZ_0 \in B\}\cap G \cap \{\exists\,\Beta^{-1}\} )}{\PP(G \cap \{\exists\,\Beta^{-1}\})}
   = \frac{\PP(\bZ_0 \in B)\PP(G \cap \{\exists\,\Beta^{-1}\})}{\PP(G \cap \{\exists\,\Beta^{-1}\})} = \PP(\bZ_0 \in B) ,
 \]
 as desired.
Note that the function \ $\varphi_\mu : \RR^d \to \CC$ \ defined in \eqref{help1} is nothing else but the characteristic function of \ $\bZ_0$ \ under \ $\PP$ \ (or \ $\QQ$).

{\sl Step 2:}
Next, we show that for each \ $r \in \ZZ_+$, \ we have
 \begin{equation}\label{sum_Delta_Y}
  \sum_{j=0}^r \bP^j \bB_{n-j} \Delta \bU_{n-j}  \to \sum_{j=0}^r \bP^j \bZ_j \qquad
  \text{$\cF_\infty$-mixing under \ $\QQ$ \ as \ $n \to \infty$.}
 \end{equation}
Let \ $r\in\ZZ_+$ \ be fixed in this step.
Since \ $\sum_{j=0}^r \bP^j \bZ_j$ \ and \ $\cF_\infty$ \ are independent under \ $\QQ$, \ we need to check that
 \begin{equation}\label{sum_Delta_Y_stab}
  \sum_{j=0}^r \bP^j \bB_{n-j} \Delta \bU_{n-j}  \to \sum_{j=0}^r \bP^j \bZ_j \qquad
  \text{$\cF_\infty$-stably under \ $\QQ$ \ as \ $n \to \infty$,}
 \end{equation}
 see the discussion after Definition 3.15 in H\"ausler and Luschgy \cite{HauLus} (or Definition \ref{Def_HL_stable_conv}).
For this, by Corollary 3.19 in H\"ausler and Luschgy \cite{HauLus} (see, also Theorem \ref{Thm_HL_Cor3_19})
 with \ $\cG := \cF_\infty$ \ and \ $\cE := \bigcup_{n\in\ZZ_+} \cF_n$, \ it is enough to show that
 \begin{equation}\label{Corollary_3_19}
  \int_\Omega
   \bbone_F
   \exp\biggl\{\ii\biggl\langle\btheta,
                                 \sum_{j=0}^r \bP^j \bB_{n-j}
                                  \Delta\bU_{n-j}\biggr\rangle\biggr\}\,
     \dd\QQ
  \to \int_\Omega
         \bbone_F
                 \exp\left\{\ii\biggl\langle\btheta, \sum_{j=0}^r \bP^j \bZ_j\biggr\rangle\right\}
                  \dd\QQ
 \end{equation}
 as \ $n \to \infty$ \ for all \ $\btheta \in \RR^d$ \ and \ $F \in \cE$.
\ Indeed, \ $\cE\subset \cF_\infty$, \ $\cE$ \ is closed under finite intersections, \ $\Omega\in \cE$ \ and \ $\sigma(\cE)=\cF_\infty$.
\ Now we turn to prove \eqref{Corollary_3_19}.
For all \ $\btheta\in\RR^d$ \ and \ $F\in\cE$, \ we have
 \[
   \exp\biggl\{\ii\biggl\langle\btheta,
                               \sum_{j=0}^r \bP^j \bB_{n-j} \Delta\bU_{n-j}\biggr\rangle\biggr\}
   = \prod_{j=0}^r \ee^{\ii\langle\btheta, \bP^j \bB_{n-j} \Delta\bU_{n-j}\rangle}
 \]
 and
 \begin{align*}
   \int_\Omega
         \bbone_F
                 \exp\left\{\ii\biggl\langle\btheta, \sum_{j=0}^r \bP^j \bZ_j\biggr\rangle\right\}
                  \dd\QQ
  & = \QQ(F) \prod_{j=0}^r \EE_\QQ(\ee^{\ii\langle\btheta,\bP^j\bZ_j\rangle})
    = \QQ(F) \prod_{j=0}^r \varphi_\mu((\bP^\top)^j \btheta) \\
  & = \int_F \prod_{j=0}^r \varphi_\mu((\bP^\top)^j \btheta) \,\dd \QQ,
 \end{align*}
 where we used that \ $\bZ_j$, \ $j\in\ZZ_+$, \ and \ $\cF_\infty$ \ are independent under \ $\QQ$,
 \ $\bZ_j$, $j\in\ZZ_+$, are identically distributed under \ $\QQ$, \ and the notation \eqref{help1}.
Hence, fixing \ $\btheta \in \RR^d$ \ arbitrarily, and using the notation
 \ $A_{n,j}:= \exp\{\ii\langle\btheta, \bP^j \bB_{n-j} \Delta\bU_{n-j}\rangle\}$,
 \ $C_j := \varphi_\mu((\bP^\top)^j \btheta)$
 \ and \ $g_{n,r} := \prod_{j=0}^r C_j - \prod_{j=0}^r A_{n,j}$ \ for \ $n \in \NN$ \ and
 \ $j \in \{0, \ldots, r\}$, \ convergence \eqref{Corollary_3_19} means that
 \ $\int_F g_{n,r} \, \dd\QQ \to 0$ \ as \ $n \to \infty$ \ for all \ $F \in \cE$.
\ By \ $|g_{n,r}| \leq 2$ \ and \eqref{MCT1}, we get
 \[
   \biggl|\int_F g_{n,r} \, \dd \QQ - \int_F L_{n-r-1} g_{n,r} \, \dd\PP\biggr|
   \leq 2 \int_F \biggl|\frac{\dd\QQ}{\dd\PP} - L_{n-r-1}\biggr| \dd\PP
   \to 0
 \]
 as \ $n \to \infty$.
\ Consequently, in order to show \eqref{Corollary_3_19}, it is enough to
 verify that \ $\lim_{n\to\infty} \int_F L_{n-r-1} g_{n,r} \, \dd\PP = 0$.
\ The condition \ $F \in \cE$ \ yields the existence of \ $n_0 \in \ZZ_+$ \ such that \ $F \in \cF_{n_0}$, \
 and consequently \ $F\in \cF_n$ \ for \ $n\geq n_0$.
\ For each \ $n \in \NN$ \ and \ $j \in \{0, \ldots, r\}$, \ put
 \[
    D_{n,j}
    := \begin{cases}
        \prod_{k=1}^r A_{n,k} & \text{if \ $j = 0$,} \\[1mm]
        \bigl(\prod_{k=0}^{j-1} C_k\bigr) \bigl(\prod_{k=j+1}^r A_{n,k}\bigr)
         & \text{if \ $1 \leq j \leq r-1$,} \\[1mm]
        \prod_{k=0}^{r-1} C_k & \text{if \ $j = r$.}
       \end{cases}
 \]
Then for each \ $n \in \NN$, \ we have
 \begin{align*}
  g_{n,r} &= \prod_{k=0}^r C_k - \prod_{k=0}^r A_{n,k} \\
      &= \prod_{k=0}^r C_k - \Biggl(\prod_{k=0}^{r-1} C_k\Biggr) A_{n,r}
         + \sum_{j=1}^{r-1}
            \Biggl[\Biggl(\prod_{k=0}^j C_k\Biggl) \Biggr(\prod_{k=j+1}^r A_{n,k} \Biggr)
                   - \Biggl(\prod_{k=0}^{j-1} C_k\Biggl)
                     \Biggr(\prod_{k=j}^r A_{n,k}\Biggr)\Biggr] \\
      &\quad + C_0 \biggr(\prod_{k=1}^r A_{n,k}\biggr) - \prod_{k=0}^r A_{n,k} \\
      &= \sum_{j=0}^r D_{n,j} (C_j - A_{n,j}) ,
 \end{align*}
 see also Lemma 8.4 in H\"ausler and Luschgy \cite{HauLus}.
Moreover, for each \ $n \in \NN$ \ and \ $j \in \{0, \ldots, r\}$, \ we have \ $|D_{n,j}| \leq 1$,
 \ and the \ $\cF_{n-j}$-measurability of \ $A_{n,j}$ \ yields that \ $D_{n,j}$ \ is \ $\cF_{n-j-1}$-measurable.
Further, for each \ $n \geq n_0 + r + 1$, \ the random variable \ $\bbone_F L_{n-r-1}$
 \ is \ $\cF_{n-r-1}$-measurable, and hence \ $\cF_{n-j-1}$-measurable for each \ $j \in \{0, \ldots, r\}$.
\ Indeed, since \ $n-r-1 \geq n_0$ \ and \ $F\in\cF_{n_0}$, \ we have \ $F\in \cF_{n-r-1}$, \ i.e.,
 \ $\bbone_F$ \ is \ $\cF_{n-r-1}$-measurable, so the \ $\cF_{n-r-1}$-measurability of \ $L_{n-r-1}$ \ yields that
 \ $\bbone_F L_{n-r-1}$ \ is \ $\cF_{n-r-1}$-measurable.
By the definition of conditional expectation, for each \ $n \geq n_0 + r + 1$, \ we obtain
 \begin{align*}
  \biggl|\int_F L_{n-r-1} g_{n,r} \, \dd\PP\biggr|
  &= \biggl|\sum_{j=0}^r \int_F L_{n-r-1} D_{n,j} (C_j - A_{n,j}) \, \dd\PP\biggr|
   = \biggl|\sum_{j=0}^r \EE_\PP\big( \bbone_F L_{n-r-1} D_{n,j} (C_j - A_{n,j}) \big) \biggr| \\
  & = \biggl|\sum_{j=0}^r \EE_\PP\big( \EE_\PP (\bbone_F L_{n-r-1} D_{n,j} (C_j - A_{n,j}) \mid \cF_{n-j-1}) \big) \biggr| \\
  & = \biggl|\sum_{j=0}^r \EE_\PP\big( \bbone_F L_{n-r-1} D_{n,j} (C_j - \EE_\PP( A_{n,j} \mid \cF_{n-j-1}) ) \big) \biggr|
 \end{align*}
 \begin{align*}
  &= \biggl|\sum_{j=0}^r
             \int_F
              L_{n-r-1} D_{n,j} (C_j - \EE_{\PP}( A_{n,j} \mid \cF_{n-j-1}))
              \, \dd\PP\biggr| .
 \end{align*}
Since \ $L_n \leq 1/\PP(G\cap\{\exists\,\Beta^{-1}\})$, \ $|C_j | \leq 1$, \ $|A_{n,j}|\leq 1$, \
 and \ $\frac{\dd\QQ}{\dd \PP} = \bbone_{G\cap\{\exists\,\Beta^{-1}\}}/\PP(G\cap\{\exists\,\Beta^{-1}\})$ \ (see the second equality in \eqref{MCT1})
 for each \ $n\geq n_0 + r +1$, \ we have
 \begin{align*}
 &\biggl|\int_F L_{n-r-1} g_{n,r} \, \dd\PP\biggr|
  \leq \sum_{j=0}^r
         \int_\Omega L_{n-r-1} |C_j - \EE_{\PP}( A_{n,j} \mid \cF_{n-j-1})| \, \dd\PP \\
  &\quad \leq \sum_{j=0}^r
         \int_{G\cap\{\exists\,\Beta^{-1}\}} \frac{1}{\PP(G\cap\{\exists\,\Beta^{-1}\})} |C_j - \EE_{\PP}( A_{n,j} \mid \cF_{n-j-1})| \, \dd\PP \\
  &\quad\phantom{\leq} + \sum_{j=0}^r \int_{\Omega\setminus (G\cap\{\exists\,\Beta^{-1}\})} L_{n-r-1}\big(|C_j| + \EE_{\PP}( |A_{n,j}| \mid \cF_{n-j-1})\big) \, \dd\PP\\
  &\quad\leq \sum_{j=0}^r \int_\Omega |C_j - \EE_{\PP}( A_{n,j} \mid \cF_{n-j-1})| \, \dd\QQ
        + 2 \sum_{j=0}^r \int_{\Omega\setminus(G\cap\{\exists\,\Beta^{-1}\})} L_{n-r-1} \, \dd\PP .
 \end{align*}
For each \ $j \in \{0, \ldots, r\}$, \ condition (iv) yields
 \begin{equation}\label{conv_C-B}
  \int_\Omega |C_j - \EE_{\PP}( A_{n,j} \mid \cF_{n-j-1})| \, \dd\QQ \to 0 \qquad
  \text{as \ $n \to \infty$.}
 \end{equation}
Indeed, since \ $|\EE_{\PP}( A_{n,j} \mid \cF_{n-j-1})|\leq 1$, \
 the family \ $\{\EE_{\PP}( A_{n,j} \mid \cF_{n-j-1}) : n \in \NN\}$ \ is uniformly integrable under \ $\QQ$ \
 for each \ $j\in\{0,\ldots,r\}$, \ and, by (iv),
 \ $\EE_{\PP}( A_{n,j} \mid \cF_{n-j-1}) \stochQ C_j$ \ as \ $n \to \infty$ \ for each \ $j \in \{0, \ldots, r\}$,
 \ so the momentum convergence theorem yields \eqref{conv_C-B}.
Further, using \eqref{MCT2} and that \ $0\leq L_{n-r-1}\leq 1/\PP(G\cap\{\exists\,\Beta^{-1}\})$, \ the dominated convergence theorem  yields that
 \[
   \int_{\Omega\setminus(G\cap\{\exists\,\Beta^{-1}\})} L_{n-r-1} \, \dd\PP
   \to
   \int_{\Omega\setminus(G\cap\{\exists\,\Beta^{-1}\})} \frac{\bbone_{G\cap\{\exists\,\Beta^{-1}\}}}{\PP(G\cap\{\exists\,\Beta^{-1}\})} \, \dd\PP
   = \QQ(\Omega \setminus (G \cap \{\exists\,\Beta^{-1}\})) = 0
 \]
 as \ $n \to \infty$.
 \ Consequently, we conclude \ $\lim_{n\to\infty} \int_F L_{n-r-1} g_{n,r} \, \dd\PP = 0$ \ for all \ $F\in\cE$, \
  and hence \eqref{Corollary_3_19}, which, as it was explained, implies \eqref{sum_Delta_Y}.

{\sl Step 3:}
Next, we check that for each \ $r\in\ZZ_+$,
 \begin{equation}\label{sum_Delta_Y_2}
  \bB_n (\bU_n - \bU_{n-r-1}) \to \sum_{j=0}^r \bP^j \bZ_j \qquad
  \text{$\cF_\infty$-mixing under \ $\QQ$ \ as \ $n \to \infty$.}
 \end{equation}
For each \ $r\in\ZZ_+$ \ and \ $j \in \{0, \ldots, r\}$, \ we have
  \begin{align}\label{help2}
   \bP^j \bB_{n-j} \Delta \bU_{n-j} - \bB_n \Delta \bU_{n-j}
    = (\bP^j - \bB_n \bB_{n-j}^{-1}) \bB_{n-j} \Delta \bU_{n-j}
    \stochQ \bzero \qquad \text{as \ $n\to\infty$.}
 \end{align}
Indeed, \ $\bB_n$ \ is invertible for sufficiently large $n\in\NN$, and
 \ $\bP^j - \bB_n \bB_{n-j}^{-1} \stochQ \bzero$ \ as \ $n \to \infty$, \ since
 for all \ $\vare>0$, \ by condition (iii),
 \begin{align*}
  \QQ(\Vert \bP^j - \bB_n \bB_{n-j}^{-1} \Vert > \vare )
   & = \frac{\PP(\{\Vert \bP^j - \bB_n \bB_{n-j}^{-1} \Vert > \vare \}\cap G\cap\{\exists\,\Beta^{-1}\} )}{\PP(G\cap\{\exists\,\Beta^{-1}\})} \\
   & \leq \frac{\PP( \{\Vert \bP^j - \bB_n \bB_{n-j}^{-1} \Vert > \vare\} \cap G) }{\PP(G\cap\{\exists\,\Beta^{-1}\})}\\
   & = \PP_G(\Vert \bP^j - \bB_n \bB_{n-j}^{-1} \Vert > \vare) \frac{\PP(G)}{\PP(G\cap \{\exists\,\Beta^{-1}\})}
    \to 0 \qquad \text{as \ $n\to\infty$.}
 \end{align*}
Further, by \eqref{sum_Delta_Y} with \ $r = 0$ \ and using the fact that \ $\cF_\infty$-mixing convergence under \ $\QQ$ \ yields convergence in
 distribution under \ $\QQ$, \ we have \ $\bB_n \Delta \bU_n \distrQ \bZ_0$ \ as \ $n \to \infty$, \
 and especially, for each \ $j\in\{0,\ldots,r\}$, \ $\bB_{n-j} \Delta \bU_{n-j} \distrQ \bZ_0$ \ as \ $n \to \infty$.
\ By Slutsky's lemma, we have \eqref{help2}.
Hence for each \ $r\in\ZZ_+$, \ we have
 \[
   \sum_{j=0}^r \bP^j \bB_{n-j}\Delta \bU_{n-j}
     - \sum_{j=0}^r \bB_n \Delta \bU_{n-j}
     \stochQ \bzero
    \qquad \text{as \ $n\to\infty$.}
 \]
Consequently, since \ $\sum_{j=0}^r \bB_n \Delta \bU_{n-j} = \bB_n (\bU_n - \bU_{n-r-1})$, \ $n \in \NN$,
 \ by \eqref{sum_Delta_Y_stab} and part (a) of Theorem 3.18 in H\"ausler and Luschgy \cite{HauLus} (see also Theorem \ref{Thm_HL_Thm3_18}),
 for each \ $r\in\ZZ_+$, \ we have
 \[
    \bB_n (\bU_n  - \bU_{n-r-1}) \to \sum_{j=0}^r \bP^j \bZ_j \qquad
      \text{$\cF_\infty$-stably under \ $\QQ$ \ as \ $n \to \infty$.}
 \]
Since \ $\sum_{j=0}^r \bP^j \bZ_j$ \ and \ $\cF_\infty$ \ are independent under \ $\QQ$ \
 (following from the \ $\QQ$-independence of \ $Z_j$, \ $j\in\ZZ_+$, \ and \ $\cF_\infty$, \
 which was proved in Step 1), by the discussion after Definition 3.15 in H\"ausler and Luschgy \cite{HauLus}
 (see also Definition \ref{Def_HL_stable_conv}), we have \eqref{sum_Delta_Y_2}.

{\sl Step 4:}
Now we turn to prove \eqref{conv_BU}.
Lemma \ref{series_convergence}, the invertibility of \ $\bP$, \ $\varrho(\bP)<1$,
 \ the condition \ $\int_{\RR^d} \log^+(\Vert\bx\Vert)\,\mu(\dd\bx)<\infty$ \
 and the fact that \ $\QQ$ \ is absolutely continuous with respect to \ $\PP$ \ (see Step 1) yield the \ $\PP$-almost sure
 and the \ $\QQ$-almost sure absolute convergence of the series \ $\sum_{j=0}^\infty \bP^j \bZ_j$.
\ Especially,
 \begin{align*}
  \sum_{j=0}^r \bP^j \bZ_j \to \sum_{j=0}^\infty \bP^j \bZ_j \qquad
  \text{as \ $r\to\infty$ \ $\QQ$-almost surely,}
 \end{align*}
 and hence
  \begin{align*}
  \sum_{j=0}^r \bP^j \bZ_j \distrQ \sum_{j=0}^\infty \bP^j \bZ_j \qquad
  \text{as \ $r\to\infty$.}
 \end{align*}
Consequently, using that \ $\sum_{j=0}^r \bP^j \bZ_j$ \ and \ $\cF_\infty$ \ are independent under \ $\QQ$ \ for each \ $r\in\ZZ_+$,
 \ by Exercise 3.4 in H\"ausler and Luschgy \cite{HauLus}, we have
  \begin{align}\label{help3}
  \sum_{j=0}^r \bP^j \bZ_j \to \sum_{j=0}^\infty \bP^j \bZ_j \qquad
  \text{$\cF_\infty$-mixing under \ $\QQ$ \ as \ $r \to \infty$.}
 \end{align}
Since \ $\bB_n \bU_n - \bB_n (\bU_n - \bU_{n-r-1}) = \bB_n \bU_{n-r-1}$, \
 and \ $\sum_{j=0}^\infty \bP^j \bZ_j$ \ and \ $\cF_\infty$ \ are independent under \ $\QQ$
 \ (following from the fact that \ $\bZ_j$, \ $j\in\ZZ_+$, \ and \ $\cF_\infty$ \ are independent under \ $\QQ$, \ which we checked in Step 1),
 \ by \eqref{sum_Delta_Y_2}, \eqref{help3} and Theorem 3.21 in H\"ausler and Luschgy \cite{HauLus} (see, also Theorem \ref{Thm_HL_Thm3_21}),
 we obtain \eqref{conv_BU} if we can check
 \begin{equation}\label{Q_n_Y_n-j-1}
  \lim_{r\to\infty} \limsup_{n\to\infty} \QQ(\|\bB_n \bU_{n-r-1}\| > \vare) = 0
 \end{equation}
 for all \ $\vare \in (0, \infty)$.
\ Since \ $\bB_n$ \ and \ $\bQ_n$ \ are invertible for sufficiently large \ $n \in \NN$, \ and \ $\bP$ \ is invertible,
 for each \ $r\in\ZZ_+$ \ and for sufficiently large \ $n\in\NN$, \ we have
 \[
   \|\bB_n \bU_{n-r-1}\|
   \leq \|\bP^{r+1}\| \cdot \|\bP^{-r-1} \bB_n \bB_{n-r-1}^{-1}\| \cdot \|\bB_{n-r-1} \bQ_{n-r-1}^{-1}\|
        \cdot \|\bQ_{n-r-1} \bU_{n-r-1}\| .
 \]
Since for each \ $r\in\ZZ_+$, \ $\bB_n \bB_{n-r-1}^{-1} \stochQ \bP^{r+1}$ \ as \ $n\to\infty$ \ (see Step 3), and
 \[
 \Vert \bP^{-r-1} \bB_n \bB_{n-r-1}^{-1} - \bI_d \Vert\leq \Vert \bP^{-r-1}\Vert \Vert \bB_n \bB_{n-r-1}^{-1} - \bP^{r+1} \Vert,
 \]
 we have \ $\bP^{-r-1} \bB_n \bB_{n-r-1}^{-1} \stochQ \bI_d$ \ as \ $n\to\infty$ \ for each \ $r\in\ZZ_+$.
\ Hence for all \ $\widetilde\vare >0$, \ $\kappa>0$ \ and \ $r\in\ZZ_+$, \ we have
 \begin{align}\label{help4}
  \QQ(\Vert \bP^{-r-1} \bB_n \bB_{n-r-1}^{-1} - \bI_d\Vert \geq  \widetilde\vare)<\kappa
   \qquad \text{for sufficiently large \ $n\in\NN$.}
 \end{align}
Consequently, with the notation \ $G_{n,r,\widetilde\vare}:=\{\Vert \bP^{-r-1} \bB_n \bB_{n-r-1}^{-1} - \bI_d\Vert < \widetilde\vare\}$, \
for all \ $\vare, \widetilde\vare, \delta,\kappa \in (0, \infty)$, \ $r\in\ZZ_+$, \ and for sufficiently large \ $n\in \NN$, \ we have
 \begin{align*}
  &\QQ(\|\bB_n \bU_{n-r-1}\| > \vare)\\
  & \leq \QQ\biggl(\|\bP^{r+1}\| \cdot \|\bP^{-r-1}\bB_n\bB_{n-r-1}^{-1}\| \cdot \| \bB_{n-r-1} \bQ_{n-r-1}^{-1} \| \cdot \|\bQ_{n-r-1} \bU_{n-r-1}\| > \vare\biggr) \\
  & = \QQ\biggl( \biggl\{ \|\bP^{r+1}\| \cdot \|\bP^{-r-1}\bB_n\bB_{n-r-1}^{-1}\| \cdot \| \bB_{n-r-1} \bQ_{n-r-1}^{-1} \| \cdot \|\bQ_{n-r-1} \bU_{n-r-1}\| > \vare \biggr\}
                   \cap G_{n,r,\widetilde\vare}  \biggr)\\
  &\phantom{=\;} + \QQ\biggl(
                    \biggl\{\|\bP^{r+1}\| \cdot \|\bP^{-r-1}\bB_n\bB_{n-r-1}^{-1}\| \cdot \| \bB_{n-r-1} \bQ_{n-r-1}^{-1} \| \cdot \|\bQ_{n-r-1} \bU_{n-r-1}\| > \vare \biggr\}
                    \cap G_{n,r,\widetilde\vare}^c  \biggr) \\
  & \leq \QQ\biggl(\!\biggl\{\|\bP^{r+1}\|\! \cdot \! \|\bP^{-r-1}\bB_n\bB_{n-r-1}^{-1} - \bI_d\| \!\cdot\! \| \bB_{n-r-1} \bQ_{n-r-1}^{-1} \|
                       \!\cdot\! \|\bQ_{n-r-1} \bU_{n-r-1}\| > \frac{\vare}{2}\biggr\}
                      \cap G_{n,r,\widetilde\vare} \biggr)\\
  &\phantom{\leq} + \QQ\biggl(\biggl\{  \|\bP^{r+1}\| \cdot \| \bB_{n-r-1} \bQ_{n-r-1}^{-1} \| \cdot \|\bQ_{n-r-1} \bU_{n-r-1}\| > \frac{\vare}{2} \biggr\}
                        \cap G_{n,r,\widetilde\vare} \biggr)\\
  &\phantom{\leq} +\QQ\biggl(\biggl\{ \|\bP^{r+1}\| \cdot \|\bP^{-r-1}\bB_n\bB_{n-r-1}^{-1}\| \cdot
                             \| \bB_{n-r-1} \bQ_{n-r-1}^{-1} \| \cdot \|\bQ_{n-r-1} \bU_{n-r-1}\| > \vare\biggr\}
                      \cap  G_{n,r,\widetilde\vare}^c  \biggr)\\
  &\leq \QQ\biggl(\|\bP^{r+1}\| \cdot \| \bB_{n-r-1} \bQ_{n-r-1}^{-1} \| \cdot \|\bQ_{n-r-1} \bU_{n-r-1}\| > \frac{\vare}{2\widetilde \vare} \biggr)\\
  &\phantom{\leq} + \QQ\biggl(\|\bP^{r+1}\| \cdot \| \bB_{n-r-1} \bQ_{n-r-1}^{-1} \| \cdot \|\bQ_{n-r-1} \bU_{n-r-1}\| > \frac{\vare}{2} \biggr)\\
  &\phantom{\leq} +\QQ\biggl(\Vert \bP^{-r-1} \bB_n \bB_{n-r-1}^{-1} - \bI_d\Vert \geq \widetilde\vare \biggr),
 \end{align*}
 where \ $G_{n,r,\widetilde\vare}^c$ \ denotes the complement of \ $G_{n,r,\widetilde\vare}$.
\ Since, by \eqref{Gelfand}, \ $\|\bP^{r+1}\| \leq \bigl(\frac{1+\varrho(\bP)}{2}\bigr)^{r+1}$ for sufficiently large \ $r \in \NN$,
 \ using also \eqref{help4}, for all \ $\vare, \delta,\kappa \in (0, \infty)$, \ $\widetilde\vare\in(0,1)$, \ and for sufficiently large \ $r \in \NN$,
 \ there exists a sufficiently large \ $n(r)\in\NN$ \ (here \ $n(r)$ \ may depend on \ $\widetilde\vare$ \ and \ $\kappa$ \ as well,
 but we do not denote this dependence) such that for \ $n\geq n(r)$, \ we have
 \begin{align*}
   &\QQ(\|\bB_n \bU_{n-r-1}\| > \vare)  \\
   &\leq \QQ\biggl(\| \bB_{n-r-1} \bQ_{n-r-1}^{-1} \| \cdot \|\bQ_{n-r-1} \bU_{n-r-1}\| > \frac{\vare}{2\widetilde \vare} \biggl(\frac{2}{1+\varrho(\bP)}\biggr)^{r+1}  \biggr)\\
   &\phantom{=} + \QQ\biggl(\| \bB_{n-r-1} \bQ_{n-r-1}^{-1} \| \cdot \|\bQ_{n-r-1} \bU_{n-r-1}\| > \frac{\vare}{2} \biggl(\frac{2}{1+\varrho(\bP)}\biggr)^{r+1} \biggr)
     + \kappa\\
   &\leq 2\QQ\biggl(\| \bB_{n-r-1} \bQ_{n-r-1}^{-1} \| \cdot \|\bQ_{n-r-1} \bU_{n-r-1}\| > \frac{\vare}{2} \biggl(\frac{2}{1+\varrho(\bP)}\biggr)^{r+1} \biggr)
         + \kappa\\
   &= 2\QQ\biggl(\| \bB_{n-r-1} \bQ_{n-r-1}^{-1} \| \cdot \|\bQ_{n-r-1} \bU_{n-r-1}\| > \frac{\vare}{2} \biggl(\frac{2}{1+\varrho(\bP)}\biggr)^{r+1},
                  \| \bB_{n-r-1} \bQ_{n-r-1}^{-1} \|\leq \delta \biggr) \\
   &\phantom{=} + 2\QQ\biggl(\| \bB_{n-r-1} \bQ_{n-r-1}^{-1} \| \!\cdot \!\|\bQ_{n-r-1} \bU_{n-r-1}\| > \frac{\vare}{2} \biggl(\frac{2}{1+\varrho(\bP)}\biggr)^{r+1},
                             \| \bB_{n-r-1} \bQ_{n-r-1}^{-1} \|> \delta  \biggr)
                 +\kappa\\
  &\leq 2\QQ\biggl(\|\bQ_{n-r-1} \bU_{n-r-1}\|
            > \frac{\vare}{2\delta} \biggl(\frac{2}{1+\varrho(\bP)}\biggr)^{r+1}\biggr)
        + 2\QQ(\|\bB_{n-r-1} \bQ_{n-r-1}^{-1}\| > \delta) + \kappa.
 \end{align*}
So for all \ $\vare, \delta,\kappa \in (0, \infty)$ \ and for sufficiently large \ $r \in \NN$,
 \ there exists a sufficiently large \ $n(r)\in\NN$ \ such that for \ $n\geq n(r)$, \ we have
 \begin{align*}
  &\QQ(\|\bB_n \bU_{n-r-1}\| > \vare)  \\
  &\leq 2\sup_{\ell\in\NN} \QQ\biggl(\|\bQ_\ell \bU_\ell\|
            > \frac{\vare}{2\delta} \biggl(\frac{2}{1+\varrho(\bP)}\biggr)^{r+1}\biggr)
        + 2\QQ(\|\bB_{n-r-1} \bQ_{n-r-1}^{-1}\| > \delta, \, \|\Beta^{-1}\| \leq \delta/2)\\
   &\quad
        + 2\QQ(\|\bB_{n-r-1} \bQ_{n-r-1}^{-1}\| > \delta, \, \|\Beta^{-1}\| > \delta/2)
        + \kappa\\
  &\leq 2\sup_{\ell\in\NN} \QQ\biggl(\|\bQ_\ell \bU_\ell\|
            > \frac{\vare}{2\delta} \biggl(\frac{2}{1+\varrho(\bP)}\biggr)^{r+1}\biggr)
        + 2\QQ\bigl(\bigl\||\bB_{n-r-1} \bQ_{n-r-1}^{-1}\| - \|\Beta^{-1}\|\bigr| > \delta/2\bigr)\\
   &\quad
        + 2\QQ(\|\Beta^{-1}\| > \delta/2)
        + \kappa ,
 \end{align*}
 where we used that \ $\QQ(\exists\,\Beta^{-1}) = 1$.
\ Similarly as we have seen in Step 3,
 condition (i) implies \ $\bQ_n \bB_n^{-1} \stochQ \Beta$ \ as \ $n \to \infty$.
\ Indeed, since \ $\PP(G)>0$, \ for all \ $\gamma>0$, \ we have
 \begin{align*}
  \QQ(\Vert \bQ_n \bB_n^{-1} - \Beta \Vert > \gamma )
   & = \frac{\PP(\{\Vert \bQ_n \bB_n^{-1} - \Beta  \Vert > \gamma \}\cap G\cap\{\exists\,\Beta^{-1}\} )}{\PP(G\cap\{\exists\,\Beta^{-1}\})} \\
   & \leq \frac{\PP( \{\Vert \bQ_n \bB_n^{-1} - \Beta \Vert > \gamma\} \cap G) }{\PP(G\cap\{\exists\,\Beta^{-1}\})}\\
   & = \PP_G(\Vert \bQ_n \bB_n^{-1} - \Beta \Vert > \gamma) \frac{\PP(G)}{\PP(G\cap \{\exists\,\Beta^{-1}\})}
    \to 0 \qquad \text{as \ $n\to\infty$.}
 \end{align*}
Since \ $\bQ_n$ \ is invertible for sufficiently large \ $n\in\NN$,
 \ $\QQ(\exists\,\Beta^{-1}) = 1$ \ and the norm function is continuous,
 we get \ $\Vert \bB_n \bQ_n^{-1} \Vert \stochQ \Vert \Beta^{-1}\Vert$ \ as \ $n \to \infty$.
\ Thus, for all \ $\vare, \delta,\kappa \in (0, \infty)$ \ and for sufficiently large \ $r \in \NN$, \ we obtain
 \[
   \limsup_{n\to\infty} \QQ(\|\bB_n \bU_{n-r-1}\| > \vare)
   \! \leq \!2\sup_{\ell\in\NN} \QQ\biggl(\|\bQ_\ell \bU_\ell\|
            > \frac{\vare}{2\delta} \biggl(\frac{2}{1+\varrho(\bP)}\biggr)^{r+1}\biggr)
        + 2\QQ(\|\Beta^{-1}\| > \delta/2)
        + \kappa.
 \]
Using condition (ii) and that \ $\frac{2}{1+\varrho(\bP)}>1$, \ for all \ $\vare, \delta, \kappa \in (0, \infty)$, \ we get
 \[
   \limsup_{r\to\infty} \limsup_{n\to\infty} \QQ(\|\bB_n \bU_{n-r-1}\| > \vare) \leq 2\QQ(\|\Beta^{-1}\| > \delta/2) + \kappa.
 \]
We have \ $\QQ(\|\Beta^{-1}\| > \delta/2) \to 0$ \ as \ $\delta \to \infty$, \ hence, taking \ $\limsup_{\delta\to\infty}$ \ and
 \ $\limsup_{\kappa\downarrow 0}$, \ we obtain \eqref{Q_n_Y_n-j-1} for all \ $\vare \in (0, \infty)$, \ and then we conclude \eqref{conv_BU}.

{\sl Step 5:} Now we turn to prove \eqref{conv_QU}.
As we have seen in Step 4, condition (i) implies \ $\bQ_n \bB_n^{-1} \stochQ \Beta$ \ as \ $n \to \infty$.
\ Hence, since \ $\Beta$ \ is \ $\cF_\infty$-measurable, by \eqref{conv_BU} (which was proved in Step 4) and parts (b) and (c) of
 Theorem 3.18 in H\"ausler and Luschgy \cite{HauLus} (see, also Theorem \ref{Thm_HL_Thm3_18}), we have
 \[
   \bQ_n \bU_n = (\bQ_n \bB_n^{-1}) (\bB_n \bU_n) \to \Beta \sum_{j=0}^\infty \bP^j \bZ_j \qquad
   \text{$\cF_\infty$-stably under \ $\QQ=\PP_{G \cap \{\exists\,\Beta^{-1}\}}$ \ as \ $n \to \infty$,}
 \]
 yielding \eqref{conv_QU}.
\proofend

\medskip

\noindent{\bf Proof of Corollary \ref{Cor_Gauss}.}
First, note that \ $\log^+(\Vert \bx\Vert) \leq \Vert \bx\Vert$, \ $\bx\in\RR^d$, \ so
 \[
   \int_{\RR^d} \log^+(\Vert \bx\Vert) \,\mu(\dd \bx) \leq \int_{\RR^d} \Vert \bx\Vert \,\mu(\dd \bx)<\infty,
 \]
 and then we can indeed apply Theorem \ref{MSLTES} and \ $\EE_\PP(\log^+(\Vert \bZ_0\Vert))<\infty$.
\ It remains to check that \ $\sum_{j=0}^\infty \bP^j \bZ_j$ \ is a \ $d$-dimensional normally distributed random variable
 with mean vector \ $\bzero\in\RR^d$ \ and covariance matrix \ $\sum_{j=0}^\infty \bP^j \bD (\bP^j)^\top$.
\ Since \ $\bP$ \ is invertible, \ $\varrho(\bP)<1$ \ and \ $\EE_\PP(\log^+(\Vert \bZ_0\Vert))<\infty$,
 by Lemma \ref{series_convergence}, we have that the series \ $\sum_{j=0}^\infty \bP^j \bZ_j$ \ is absolutely convergent \ $\PP$-a.s.,
 and hence, by the continuity theorem, we get
 \begin{align*}
  \EE_\PP\big(\ee^{\ii \langle \btheta, \sum_{j=0}^\infty \bP^j \bZ_j\rangle}\big)
   & = \lim_{r\to\infty} \EE_\PP\big(\ee^{\ii \langle \btheta, \sum_{j=0}^r \bP^j \bZ_j\rangle}\big)
     = \lim_{r\to\infty} \prod_{j=0}^r \EE_\PP\big(\ee^{\ii \langle (\bP^j)^\top \btheta, \bZ_j\rangle}\big)\\
   & = \lim_{r\to\infty} \prod_{j=0}^r \ee^{-\frac{1}{2} \langle \bD (\bP^j)^\top \btheta , (\bP^j)^\top \btheta \rangle}
     = \ee^{-\frac{1}{2} \left\langle  \left( \sum_{j=0}^\infty \bP^j \bD (\bP^j)^\top \right) \btheta, \btheta \right\rangle},
      \qquad \btheta\in\RR^d,
 \end{align*}
 where the series \ $\sum_{j=0}^\infty \bP^j \bD (\bP^j)^\top$ \ is absolutely convergent, since,
 by \eqref{Gelfand},
 \[
  \sum_{j=0}^\infty \Vert \bP^j \bD (\bP^j)^\top \Vert
    \leq \sum_{j=0}^\infty  \Vert \bP^j \Vert \Vert \bD \Vert \Vert (\bP^j)^\top \Vert
    \leq \Vert \bD \Vert \sum_{j=0}^{k_0-1}  \Vert \bP^j \Vert^2
          + \Vert \bD \Vert \sum_{j=k_0}^\infty  \left(\frac{1+\varrho(\bP)}{2}\right)^{2j}
    <\infty,
 \]
 where \ $k_0$ \ is appearing in \eqref{Gelfand}.
So \ $\sum_{j=0}^\infty \bP^j \bZ_j$ \ is a \ $d$-dimensional normally distributed random variable
 with mean vector \ $\bzero\in\RR^d$ \ and covariance matrix \ $\sum_{j=0}^\infty \bP^j \bD (\bP^j)^\top$,
 \ as desired.
\proofend

\medskip

\noindent{\bf Proof of Corollary \ref{Cor_Stable}.}
First, note that the integral appearing in \eqref{help_Stable_limit_char} is convergent, since, by Cauchy-Schwarz's inequality and \eqref{Gelfand},
 for all \ $\btheta\in\RR^d$,
 \begin{align*}
  & \int_{S_{d-1}}\sum_{j=0}^\infty \big\vert \langle (\bP^j)^\top \btheta,\bx \rangle \big\vert^\alpha \, \Pi(\dd \bx)
    \leq \int_{S_{d-1}}\sum_{j=0}^\infty \big\Vert (\bP^j)^\top \btheta \Vert^\alpha \Vert \bx\Vert^\alpha\, \Pi(\dd \bx) \\
  &\leq \big\Vert \btheta \Vert^\alpha \int_{S_{d-1}} \sum_{j=0}^\infty \Vert \bP^j\Vert^\alpha \Vert \bx\Vert^\alpha \, \Pi(\dd \bx) \\
  &\leq  \big\Vert \btheta \Vert^\alpha \int_{S_{d-1}} \sum_{j=0}^{k_0-1} \Vert \bP^j\Vert^\alpha \Vert \bx\Vert^\alpha \, \Pi(\dd \bx)
           + \big\Vert \btheta \Vert^\alpha \int_{S_{d-1}} \sum_{j=k_0}^\infty \left(\frac{1+\varrho(\bP)}{2}\right)^{\alpha j} \Vert \bx\Vert^\alpha \, \Pi(\dd \bx)\\
  &= \big\Vert \btheta \Vert^\alpha \sum_{j=0}^{k_0-1} \Vert \bP^j\Vert^\alpha \, \Pi(S_{d-1})
      + \big\Vert \btheta \Vert^\alpha \sum_{j=k_0}^\infty \left(\frac{1+\varrho(\bP)}{2}\right)^{\alpha j} \, \Pi(S_{d-1})
    <\infty,
 \end{align*}
 where \ $k_0$ \ is appearing in \eqref{Gelfand} and we also used that \ $\frac{1+\varrho(\bP)}{2}\in(0,1)$ \ and \ $\Pi(S_{d-1})<\infty$.

Next, we check that \ $\EE_\PP(\log^+(\Vert \bzeta\Vert))<\infty$.
\ We have
 \begin{align*}
   \EE_\PP(\log^+(\Vert \bzeta\Vert))
    & = \EE_\PP(\log(\Vert \bzeta\Vert)\bbone_{\{ \Vert \bzeta\Vert \geq 1\}} )
      = \int_0^\infty \PP( \log(\Vert \bzeta\Vert)\bbone_{\{ \Vert \bzeta\Vert \geq 1\}}  \geq y)\,\dd y\\
    & = \int_0^\infty \PP( \log(\Vert \bzeta\Vert) \geq y, \Vert \bzeta\Vert \geq 1)\,\dd y
        = \int_0^\infty \PP( \Vert \bzeta\Vert \geq \ee^y, \Vert \bzeta\Vert \geq 1 ) \,\dd y \\
    & = \int_0^1 \PP( \Vert \bzeta\Vert \geq \ee^y) \,\dd y
        + \int_1^\infty \PP( \Vert \bzeta\Vert \geq \ee^y) \,\dd y  \\
    & \leq 1 + \int_\ee^\infty \PP( \Vert \bzeta\Vert \geq z) \frac{1}{z}\,\dd z.
 \end{align*}
Since \ $\bzeta$ \ has a \ $d$-dimensional stable distribution, it belongs to its own domain of attraction,
 and then it is known that the function \ $\RR_{++}\ni z \mapsto \PP(\Vert \bzeta\Vert \geq z)$ \ is regularly varying with tail index
 \ $\alpha$.
\ As a consequence, the function \ $\RR_{++}\ni z \mapsto z^\alpha\PP(\Vert \bzeta\Vert \geq z)=:L(z)$ \ is slowly varying.
Hence there exists \ $z_0\in(\ee,\infty)$ \ such that \ $z^{-\frac{\alpha}{2}} L(z)\leq 1$ \ for all \ $z\in[z_0,\infty)$, \ see, e.g.,
 Bingham et al.\ \cite[Proposition 1.3.6.(v)]{BinGolTeu}.
Consequently, we have
 \begin{align*}
   &\int_\ee^\infty \PP( \Vert \bzeta\Vert \geq z) \frac{1}{z}\,\dd z
      = \int_\ee^{z_0} z^{-\alpha} L(z) \frac{1}{z}\,\dd z
         + \int_{z_0}^\infty z^{-\alpha} L(z) \frac{1}{z}\,\dd z\\
  &\leq \int_\ee^{z_0} z^{-\alpha} L(z) \frac{1}{z}\,\dd z
         + \int_{z_0}^\infty z^{-\frac{\alpha}{2}} \frac{1}{z}\,\dd z
   \leq \int_\ee^{z_0} \frac{1}{z}\,\dd z
         + \int_{z_0}^\infty z^{-1-\frac{\alpha}{2}} \,\dd z
   <\infty,
 \end{align*}
 since \ $z^{-\alpha} L(z) = \PP(\Vert \bzeta\Vert > z)\leq 1$, \ $z\in\RR_{++}$.

Hence one can indeed apply Theorem \ref{MSLTES} and \ $\EE_\PP(\log^+(\Vert \bZ_0\Vert))<\infty$.
\ It remains to check that the characteristic function of \ $\sum_{j=0}^\infty \bP^j \bZ_j$ \ is given by \eqref{help_Stable_limit_char}.
Since \ $\bP$ \ is invertible, $\varrho(\bP)<1$, \ and \ $\EE_\PP(\log^+(\Vert \bZ_0\Vert))<\infty$, \ by Lemma \ref{series_convergence},
 we have that \ $\sum_{j=0}^\infty \bP^j \bZ_j$ \ is absolutely convergent \ $\PP$-a.s., and hence, by the continuity theorem, we get
 \begin{align*}
   \EE_\PP\big(\ee^{\ii \langle \btheta, \sum_{j=0}^\infty \bP^j \bZ_j\rangle}\big)
   & = \lim_{r\to\infty} \prod_{j=0}^r \EE_\PP\big(\ee^{\ii \langle (\bP^j)^\top\btheta,  \bZ_j \rangle}\big)
     = \lim_{r\to\infty} \exp\left\{ - \sum_{j=0}^r \int_{S_{d-1}} \vert \langle (\bP^j)^\top \btheta, \bx \rangle \vert^\alpha\,\Pi(\dd \bx) \right\} \\
    & = \exp\left\{ -\sum_{j=0}^\infty \int_{S_{d-1}} \vert \langle (\bP^j)^\top \btheta, \bx \rangle \vert^\alpha\,\Pi(\dd \bx)  \right\}, \qquad \btheta\in\RR^d,
 \end{align*}
 yielding \eqref{help_Stable_limit_char}.

In the special case when \ $\bzeta$ \ has a $d$-dimensional Cauchy distribution with parameter \ $(\bzero,\bI_d)$, \ we have
 \ $\int_{S_{d-1}} \vert \langle (\bP^j)^\top \btheta, \bx \rangle \vert^\alpha\,\Pi(\dd \bx) = \Vert (\bP^j)^\top \btheta \Vert$, \ $\btheta\in\RR^d$, \ $j\in\ZZ_+$,
 \ yielding \eqref{help_Cauchy_limit_char}.
\proofend

\medskip

\noindent{\bf Proof of Corollary \ref{MSLTES_discrete_S}.}
Let \ $\{\bs_k: k\in\NN\}$ \ be the range of \ $\bS$, \ let \ $G_k:=\{\bS = \bs_k\}$, \ $k\in\NN$, \ and
 \ $I:=\{ k\in\NN : \PP(G_k \cap \{\exists \; \Beta^{-1}\})>0\}$.
\ Since \ $\PP(\exists \; \Beta^{-1})>0$ \ (due to \ $G=\Omega$), \ we have that \ $I$ \ is not the empty set.
\ Further, since \ $\PP_{G_k\cap \{\exists \; \Beta^{-1}\}}$ \ is absolutely continuous with respect to \ $\PP_{\{\exists \; \Beta^{-1}\}}$, \ by (iv'),
 and using that convergence in \ $\PP_{\{\exists\, \Beta^{-1}\}}$-probability yields convergence in \ $\PP_{G_k\cap\{\exists\, \Beta^{-1}\}}$-probability
 (which can be checked similarly as in case of \ $\PP_G$ \ and \ $\PP_{G\cap \{\exists\, \Beta^{-1}\}}$ \
 as we have seen in the proof of Step 3 of Theorem \ref{MSLTES}), we have for each \ $k\in I$ \ and \ $\btheta\in\RR^d$,
 \begin{align*}
  \EE_\PP\bigl(\ee^{\ii\langle\btheta,\bB_n\Delta\bU_n\rangle}
                     \mid \cF_{n-1}\bigr)
     \stochGketa \int_{\RR^d}\ee^{\ii\langle \btheta, \bS\bx \rangle}\,\mu(\dd\bx)
                & = \int_{\RR^d}\ee^{\ii\langle \btheta, \bs_k\bx \rangle}\,\mu(\dd\bx)
                  = \EE_\PP\left( \ee^{\ii\langle \bs_k^\top \btheta,\bZ_0 \rangle} \right)\\
                & = \EE_\PP\left( \ee^{\ii\langle \btheta,\bs_k \bZ_0 \rangle} \right)
 \end{align*}
 as \ $n\to\infty$.
\ Moreover, \ $\EE_\PP(\log^+(\Vert \bs_k\bZ_0\Vert))<\infty$, \ since
 \begin{align*}
   &\log^+(\Vert \bs_k\bZ_0\Vert)
     = \log(\Vert \bs_k\bZ_0\Vert) \bbone_{\{ \Vert \bs_k\bZ_0\Vert \geq 1\}}
      \leq \log(\Vert \bs_k \Vert \Vert \bZ_0\Vert) \bbone_{\{ \Vert \bs_k \Vert \Vert \bZ_0\Vert \geq 1\}}\\
   &\quad \leq  \log(\Vert \bs_k \Vert ) \bbone_{\{\bs_k \ne \bzero\}} \bbone_{\{ \Vert \bs_k \Vert \Vert \bZ_0\Vert \geq 1\}}
            + \log(\Vert \bZ_0 \Vert)\bbone_{\{ \Vert \bZ_0\Vert \geq 1\}}
            + \log(\Vert \bZ_0 \Vert) \bbone_{\{ \frac{1}{\Vert \bs_k\Vert}\leq  \Vert \bZ_0\Vert < 1\}} \bbone_{\{\bs_k \ne \bzero\}},
 \end{align*}
 which yields that
 \[
  \EE_\PP(\log^+(\Vert \bs_k\bZ_0\Vert))
    \leq \log(\Vert \bs_k \Vert ) \bbone_{\{\bs_k \ne \bzero\}}
          + \EE_\PP(\log^+(\Vert \bZ_0\Vert)) < \infty.
 \]
Hence, by Theorem \ref{MSLTES}, for each \ $k\in I$, \ we have
 \begin{equation*}
  \bB_n \bU_n \to \sum_{j=0}^\infty \bP^j \bs_k\bZ_j = \sum_{j=0}^\infty \bP^j \bS\bZ_j\qquad
  \text{$\cF_\infty$-mixing under \ $\PP_{G_k \cap \{\exists\,\Beta^{-1}\}}$ \ as \ $n \to \infty$,}
 \end{equation*}
 and
 \begin{equation*}
   \bQ_n \bU_n \to \Beta \sum_{j=0}^\infty \bP^j \bs_k\bZ_j  = \Beta\sum_{j=0}^\infty \bP^j \bS\bZ_j \qquad
   \text{$\cF_\infty$-stably under \ $\PP_{G_k \cap \{\exists\,\Beta^{-1}\}}$ \ as \ $n \to \infty$.}
 \end{equation*}
Note that, since \ $G=\Omega$, \ we have \ $\PP(\exists\; \Beta^{-1})>0$ \ and for all \ $A\in\cF$,
 \begin{align*}
   \PP_{\{\exists\; \Beta^{-1}\}}(A)
    & = \sum_{k=1}^\infty \PP_{\{\exists\; \Beta^{-1}\}} (A\cap G_k)
     = \sum_{k=1}^\infty \frac{\PP(A\cap G_k\cap \{\exists\; \Beta^{-1}\})}{\PP(\exists\;\Beta^{-1})}\\
    & = \sum_{k\in I} \frac{\PP(A\cap G_k\cap \{\exists\; \Beta^{-1}\})}{\PP(G_k\cap\{\exists\;\Beta^{-1}\})}
                         \frac{\PP( G_k\cap\{\exists\; \Beta^{-1}\} )}{\PP(\exists\;\Beta^{-1})}
      = \sum_{k\in I} \PP_{ G_k\cap\{\exists\; \Beta^{-1}\}}(A)\PP_{\{\exists\; \Beta^{-1}\}} (G_k),
 \end{align*}
 so we have
 \[
   \PP_{\{\exists\; \Beta^{-1}\}} = \sum_{k\in I} \PP_{\{\exists\; \Beta^{-1}\}}(G_k) \PP_{G_k\cap \{\exists\; \Beta^{-1}\}},
 \]
 where \ $\sum_{k\in I} \PP_{\{\exists\; \Beta^{-1}\}}(G_k)=1$.
\ Finally, Proposition 3.24 in H\"ausler and Luschgy \cite{HauLus} (see also Theorem \ref{Thm_HL_Propostion3_24}) yields the statement.
\proofend

\appendix

\vspace*{5mm}

\noindent{\bf\Large Appendix}

\section{Stable convergence and L\'evy's upwards theorem}
\label{App_1}

First, we recall the notions of stable and mixing convergence.

\begin{Def}\label{Def_HL_stable_conv}
Let \ $(\Omega,\cF,\PP)$ \ be a probability space and \ $\cG\subset \cF$ \ be a sub-$\sigma$-field.
Let \ $(\bX_n)_{n\in\NN}$ \ and \ $\bX$ \ be \ $\RR^d$-valued random variables defined on \ $(\Omega,\cF,\PP)$,
 \ where \ $d\in\NN$.

\noindent (i) We say that \ $\bX_n$ \ converges \ $\cG$-stably to \ $\bX$ \ as \ $n\to\infty$, \ if the conditional distribution
 \ $\PP^{\bX_n\mid \cG}$ \ of \ $\bX_n$ \ given \ $\cG$ \ converges weakly to the conditional distribution
 \ $\PP^{\bX\mid \cG}$ \ of \ $\bX$ \ given \ $\cG$ \ as \ $n\to\infty$ \ in the sense of weak convergence of Markov kernels.
It equivalently means that
 \[
   \lim_{n\to\infty} \EE_\PP(\xi \EE_\PP(h(\bX_n) \mid \cG ) )
      = \EE_\PP( \xi \EE_\PP(h(\bX) \mid \cG ) )
 \]
 for all random variables \ $\xi:\Omega\to\RR$ \ with \ $\EE_\PP(\vert \xi\vert)<\infty$ \ and for all bounded and continuous functions
 \ $h:\RR^d\to\RR$.

\noindent (ii) We say that \ $\bX_n$ \ converges \ $\cG$-mixing to \ $\bX$ \ as \ $n\to\infty$, \
 if \ $\bX_n$ \ converges \ $\cG$-stably to \ $\bX$ \ as \ $n\to\infty$, \ and \ $\PP^{\bX\mid \cG} = \PP^\bX$ \ $\PP$-almost surely, where
 \ $\PP^\bX$ \ denotes the distribution of \ $\bX$ \ on \ $(\RR^d,\cB(\RR^d))$ \ under \ $\PP$.
\ Equivalently, we can say that \ $\bX_n$ \ converges \ $\cG$-mixing to \ $\bX$ \ as \ $n\to\infty$,
 \ if \ $\bX_n$ \ converges \ $\cG$-stably to \ $\bX$ \ as \ $n\to\infty$, \ and \ $\sigma(\bX)$ \ and \ $\cG$ \ are independent,
 which equivalently means that
 \[
   \lim_{n\to\infty} \EE_\PP(\xi \EE_\PP(h(\bX_n) \mid \cG ) )
      = \EE_\PP(\xi) \EE_\PP(h(\bX))
 \]
 for all random variables \ $\xi:\Omega\to\RR$ \ with \ $\EE_\PP(\vert \xi\vert)<\infty$ \ and for all bounded and continuous functions
 \ $h:\RR^d\to\RR$.
\end{Def}

In Definition \ref{Def_HL_stable_conv}, \ $\PP^{\bX_n\mid \cG}$, \ $n\in\NN$, \ and \ $\PP^{\bX\mid \cG}$ \ are
 the \ $\PP$-almost surely unique \ $\cG$-measurable Markov kernels from \ $(\Omega,\cF)$ \ to \ $(\RR^d,\cB(\RR^d))$ \ such that for each \ $n\in\NN$,
 \[
   \int_G \PP^{\bX_n\mid \cG}(\omega,B)\,\PP(\dd \omega)
     = \PP(\bX_n^{-1}(B) \cap G)
     \qquad \text{for all \ $G\in \cG$, \ $B\in\cB(\RR^d)$.}
 \]
 and
 \[
   \int_G \PP^{\bX\mid \cG}(\omega,B)\,\PP(\dd \omega)
     = \PP(\bX^{-1}(B) \cap G)
     \qquad \text{for all \ $G\in \cG$, \ $B\in\cB(\RR^d)$,}
 \]
 respectively.
For the notion of weak convergence of Markov kernels towards a Markov kernel, see H\"ausler and Luschgy \cite[Definition 2.2]{HauLus}.
For more details on stable convergence, see H\"ausler and Luschgy \cite[Chapter 3 and Appendix A]{HauLus}.
In particular, it turns out that \ $\bX_n$ \ converges \ $\cG$-stably to \ $\bX$ \ as \ $n\to\infty$
 \ if and only if \ $\lim_{n\to\infty} \EE_\PP(\xi h(\bX_n)) = \EE_\PP(\xi h(\bX))$ \ for all \ $\cG$-measurable random variables
 \ $\xi:\Omega\to\RR$ \ with \ $\EE_\PP(\vert \xi \vert)<\infty$ \ and for all bounded and continuous functions
 \ $h:\RR^d\to\RR$ \ (following from Theorem 3.17 in H\"ausler and Luschgy \cite{HauLus}).
Furthermore, \ $\bX_n$ \ converges \ $\cG$-mixing to \ $\bX$ \ as \ $n\to\infty$ \ if and only if
 \ $\lim_{n\to\infty} \EE_\PP(\xi h(\bX_n)) = \EE_\PP(\xi)\EE_\PP(h(\bX))$ \ for all \ $\cG$-measurable random variables
 \ $\xi:\Omega\to\RR$ \ with \ $\EE_\PP(\vert \xi \vert)<\infty$ \ and for all bounded and continuous functions
 \ $h:\RR^d\to\RR$ \ (following from Corollary 3.3 in H\"ausler and Luschgy \cite{HauLus}).

Next, we recall four results about stable convergence of random variables, which play important roles in the proofs of Theorem \ref{MSLTES}
 and Corollary \ref{MSLTES_discrete_S}.

\begin{Thm}[H\"ausler and Luschgy {\cite[Theorem 3.18]{HauLus}}]\label{Thm_HL_Thm3_18}
Let \ $\bX_n$, \ $n\in\NN$, \ $\bX$, \ $\bY_n$, \ $n\in\NN$, \ and \ $\bY$ \ be \ $\RR^d$-valued random variables on a probability space \ $(\Omega,\cF,\PP)$, \
  and \ $\cG\subset \cF$ \ be a sub-$\sigma$-field.
Assume that \ $\bX_n\to \bX$ \ $\cG$-stably as \ $n\to\infty$.
 \begin{itemize}
   \item[(a)] If \ $\Vert \bX_n - \bY_n\Vert \stoch 0$ \ as \ $n\to\infty$, \ then \ $\bY_n\to \bX$ \ $\cG$-stably as \ $n\to\infty$.
   \item[(b)] If \ $\bY_n \stoch \bY$ \ as \ $n\to\infty$, \ and \ $\bY$ \ is \ $\cG$-measurable, then
               \ $(\bX_n,\bY_n)\to (\bX,\bY)$ \ $\cG$-stably as \ $n\to\infty$.
   \item[(c)] If \ $g:\RR^d\to\RR^d$ \ is a Borel-measurable function such that \ $\PP^\bX(\{ \bx\in\RR^d : \text{$g$ \ is not continuous at \ $\bx$}\})=0$, \
              then \ $g(\bX_n)\to g(\bX)$ \  $\cG$-stably as \ $n\to\infty$.
              \ Here recall that \ $\PP^\bX$ \ denotes the distribution of \ $\bX$ \ on \ $(\RR^d,\cB(\RR^d))$ \ under \ $\PP$.
 \end{itemize}
\end{Thm}

\begin{Thm}[H\"ausler and Luschgy {\cite[Corollary 3.19]{HauLus}}]\label{Thm_HL_Cor3_19}
Let \ $\bX_n$, \ $n\in\NN$, \ and \ $\bX$ \ be \ $\RR^d$-valued random variables on a probability space \ $(\Omega,\cF,\PP)$, \
  and \ $\cG\subset \cF$ \ be a sub-$\sigma$-field.
Let \ $\cE\subset \cG$ \ be closed under finite intersections such that \ $\Omega\in\cE$ \ and \ $\sigma(\cE) = \cG$, \ where
 \ $\sigma(\cE)$ \ denotes the \ $\sigma$-algebra generated by \ $\cE$.
\ Then the following statements are equivalent:
 \begin{enumerate}
  \item[\textup{(i)}]
    $\bX_n\to \bX$ \ $\cG$-stably as \ $n\to\infty$,
  \item[\textup{(ii)}]
    $\lim_{n\to\infty} \EE_\PP(\bbone_F \ee^{\ii\langle \bu,\bX_n\rangle}) =  \EE_\PP(\bbone_F \ee^{\ii\langle \bu,\bX\rangle})$
    \ for all \ $F\in \cE$ \ and \ $\bu\in \RR^d$,
  \item[\textup{(iii)}]
    $\langle \bu,\bX_n\rangle \to \langle \bu,\bX\rangle$ \ $\cG$-stably for all \ $\bu\in\RR^d$.
 \end{enumerate}
\end{Thm}

\begin{Thm}[H\"ausler and Luschgy {\cite[Theorem 3.21]{HauLus}}]\label{Thm_HL_Thm3_21}
Let \ $\bX_{n,r}$, $\bX_r$, \ $n,r\in\NN$, \ $\bX$, and \ $\bY_n$, \ $n\in\NN$, \ be \ $\RR^d$-valued random variables
 on a probability space \ $(\Omega,\cF,\PP)$, \ and \ $\cG\subset \cF$ \ be a sub-$\sigma$-field.
Assume that
 \begin{enumerate}
  \item[\textup{(i)}]
    for each \ $r\in\NN$, \ we have \ $\bX_{n,r}\to \bX_r$ \ $\cG$-stably as \ $n\to\infty$,
  \item[\textup{(ii)}]
    $\bX_r\to \bX$ \ $\cG$-stably as \ $r\to\infty$,
  \item[\textup{(iii)}]
    $\lim_{r\to\infty} \limsup_{n\to\infty} \PP(\Vert \bX_{n,r} - \bY_n\Vert > \varepsilon)=0$ \ for all \ $\varepsilon>0$.
 \end{enumerate}
Then \ $\bY_n\to\bX$ \ $\cG$-stably as \ $n\to\infty$.
\end{Thm}

\begin{Thm}[H\"ausler and Luschgy {\cite[Proposition 3.24]{HauLus}}]\label{Thm_HL_Propostion3_24}
Let \ $(\Omega,\cF,\PP)$ \ be a probability space.
Suppose that \ $\PP=\sum_{i=1}^\infty p_i \QQ_i$, \ where \ $\QQ_i$, \ $i\in\NN$, \ is a probability measure on \ $(\Omega,\cF)$ \
 and \ $p_i\in[0,1]$, \ $i\in\NN$, \ satisfying \ $\sum_{i=1}^\infty p_i=1$.
\ Let \ $\bX_n$, \ $n\in\NN$, \ and \ $\bX$ \ be \ $\RR^d$-valued random variables on \ $(\Omega,\cF,\PP)$.
\ If \ $\bX_n$ \ converges \ $\cG$-stably to \ $\bX$ \ under \ $\QQ_i$ \ as \ $n\to\infty$ \ for each \ $i\in\NN$ \
 satisfying \ $p_i>0$, \ then \ $\bX_n$ \ converges \ $\cG$-stably to \ $\bX$ \ under \ $\PP$ \ as \ $n\to\infty$.
\end{Thm}

Finally, we recall L\'evy's upwards theorem used in the proof of Theorem \ref{MSLTES}.

\begin{Thm}[L\'evy's upwards theorem]\label{Thm_Levy01}
Let \ $(\Omega,\cF,\PP)$ \ be a probability space, and let \ $\xi$ \ be a real-valued random variable such that
 \ $\EE_\PP(\vert\xi\vert)<\infty$ \ and \ $(\cF_n)_{n\in\ZZ_+}$ \ be a filtration with
 \ $\cF_\infty := \sigma\big(\bigcup_{n\in\ZZ_+} \cF_n \big)$.
\ Then
 \[
   \EE_\PP(\xi\mid \cF_n) \asP \EE_\PP(\xi\mid \cF_\infty)
    \quad \text{as \ $n\to\infty$,}
    \qquad \text{and} \qquad
   \EE_\PP(\xi\mid \cF_n) \meanP \EE_\PP(\xi\mid \cF_\infty)
   \quad \text{as \ $n\to\infty$.}
 \]
\end{Thm}

We note that Theorem \ref{Thm_Levy01} sometimes is called L\'evy's zero-one law as well, since if \ $\xi=\bbone_A$, \  where \ $A\in\cF_\infty$, \
 then it yields that \ $\PP(A\mid \cF_n)\asP \bbone_A$ \ as \ $n\to\infty$, \ where the limit can be zero or one.

\section*{Acknowledgements}
We would like to thank the referee for the comments that helped us improve the paper.

\end{document}